\documentclass[hidelinks,onefignum,onetabnum]{siamart220329}

\usepackage{lipsum}
\usepackage{amsfonts}
\usepackage{graphicx}
\usepackage{epstopdf}
\usepackage{algorithmic}
\ifpdf
  \DeclareGraphicsExtensions{.eps,.pdf,.png,.jpg}
\else
  \DeclareGraphicsExtensions{.eps}
\fi

\newsiamremark{remark}{Remark}

\newsiamremark{hypothesis}{Hypothesis}
\crefname{hypothesis}{Hypothesis}{Hypotheses}
\newsiamthm{claim}{Claim}

\headers{Solving FK equation of Dean-Kawasaki models with FHT ansatz}{X. Tang, L. Ying}

\title{Solving the Fokker-Planck equation of discretized Dean-Kawasaki models with functional hierarchical tensor}

\author{
  Xun Tang\thanks{Corresponding author. Institute for Computational and Mathematical Engineering (ICME), Stanford University, Stanford, CA 94305, USA. 
  (\email{xuntang@stanford.edu})}
  \and
  Lexing Ying\thanks{Department of Mathematics and Institute for Computational and Mathematical Engineering (ICME), Stanford University, Stanford, CA 94305, USA. 
  (\email{lexing@stanford.edu})}
  \funding{X.T. and L.Y. are supported by AFOSR MURI award FA9550-24-1-0254.}
  }

\usepackage{amsopn}

\makeatletter
\newcommand*{\addFileDependency}[1]{%
  \typeout{(#1)}%
  \@addtofilelist{#1}%
  \IfFileExists{#1}{}{\typeout{No file #1.}}%
}
\makeatother

\usepackage{amsmath,amsfonts,bm}

\def\eqref#1{equation~\ref{#1}}

\def\1{\bm{1}}

\DeclareMathAlphabet{\mathsfit}{\encodingdefault}{\sfdefault}{m}{sl}
\SetMathAlphabet{\mathsfit}{bold}{\encodingdefault}{\sfdefault}{bx}{n}

\newcommand{\R}{\mathbb{R}}

\usepackage{amsmath,amssymb,amsfonts}
\usepackage{latexsym}
\usepackage{indentfirst}
\usepackage{graphicx}
\usepackage{placeins}
\usepackage{enumerate}
\usepackage{booktabs}  
\usepackage{algorithm}
\usepackage{algorithmic}
\usepackage{multirow}
\usepackage{optidef}
\usepackage{hyperref}
\usepackage{diagbox}
\usepackage{subcaption}
\usepackage[normalem]{ulem}
\newcommand{\edit}[1]{{#1}}

\usepackage{cleveref}

\begin{document}

\maketitle

\begin{abstract}
    We propose a particle-based workflow for approximating the time-dependent law of finite-volume discretizations of the Dean-Kawasaki (DK) model. After discretization, the state is a nonnegative vector whose total mass is conserved by the finite-volume scheme, and it is therefore supported on a probability simplex. To enable tensor-network density estimation, we map the simplex to an unconstrained Euclidean space using a centered logarithmic transform and then apply a wavelet transform that organizes degrees of freedom by spatial scale. On the transformed variables, we fit the probability density with a functional hierarchical tensor over a wavelet basis (FHT-W), i.e., a hierarchical-Tucker/tree-tensor-network representation of the coefficient tensor of a fixed univariate basis expansion. We illustrate the method on 1D and 2D examples with $d=64$ degrees of freedom, including cases with external potentials and pairwise interactions. The method accurately captures the site-wise correlations and other observables of the true model.
\end{abstract}



\section{Introduction}\label{sec:intro}
We study the Dean-Kawasaki model \cite{dean1996langevin,kawasaki1998microscopic} with periodic boundary conditions:
\begin{equation}\label{eqn: Dean eq}
    \partial_t \pi(x,t) = \frac{1}{\beta}\Delta \pi + \mathrm{div}\left(\sqrt{\frac{2\pi}{\beta N}}\eta + \pi\nabla V_1  + \pi\nabla V_2 * \pi\right), \quad \pi(x, 0) = \pi_0(x),
\end{equation}
where \(x \in [0,1]^n\) with \(n \in \{1,2\}\), \(\eta(x, t)\) is a \(\R^n\)-valued space-time white noise, \(N\) controls the strength of the white noise term, \(\beta\) is the inverse temperature controlling the diffusion coefficient, \(V_1(\cdot)\) is an external potential term, and \(V_2(\cdot)\) is a pairwise potential term modeling particle interactions. \Cref{eqn: Dean eq} is a stochastic partial differential equation (SPDE). In \Cref{eqn: Dean eq}, the solution term \(\pi\) models the particle density (e.g., cell density or chemical concentration density). \edit{In particular, \(\pi(x, t)\) stands for the density at location \(x\) at time \(t\).}

In this work, we study the Dean-Kawasaki model after spatial discretization, where the domain is divided into \(d\) grid cells. For example, in 1D, the discretization leads to a \(d\)-dimensional cell-average vector \(\Pi = (\Pi_1, \ldots, \Pi_d)\). Under a discretization scheme such as the finite-volume scheme to be introduced in \Cref{sec: dean}, the SPDE in \Cref{eqn: Dean eq} discretizes to a stochastic differential equation (SDE):
\begin{equation}\label{eqn: Dean SDE}
    d \Pi(t) = \bm{\mu}(\Pi(t)) dt  + \bm{\sigma}(\Pi(t)) dB_t,
\end{equation}
where \(B_t\) is the Brownian motion (i.e., Wiener process).

\Cref{eqn: Dean SDE} has an associated Fokker-Planck equation. \edit{We let a function \(P_{\pi}(\pi, t)\) denote the probability density of \(\frac{\Pi}{d}\) at time \(t\). The Fokker-Planck equation is as follows}:
\begin{equation}\label{eqn: Fokker-Planck}
    \partial_t P_{\pi}(\pi, t) = -\mathrm{div}(\bm{\mu}(\pi)P_{\pi}) + \frac{1}{2}\nabla_{\pi}^2:(\bm{\Sigma}(\pi)P_{\pi}), \quad \pi \in \Delta, t \in [0, T],
\end{equation}
where \(\bm{\Sigma}(\pi) = \bm{\sigma}(\pi)\bm{\sigma}^{\top}(\pi)\) and \(\nabla_{\pi}^2:(\bm{\Sigma}(\pi)P_{\pi}) = \sum_{ij}\partial_{\pi_i}\partial_{\pi_j}\left(P_{\pi}(\pi, t)\left(\bm{\Sigma}(\pi)\right)_{ij}\right)\), and \(\Delta\) is a simplex since entries of \(\Pi\) sums to \(d\).

In this work, our goal is to solve the high-dimensional partial differential equation (PDE) in \Cref{eqn: Fokker-Planck}. Numerically treating the discretized Dean-Kawasaki (DK) model is difficult, and we apply several techniques to simplify the task. First, we use a particle-based method whereby one solves the PDE by performing density estimation. After simulating the SDE in \Cref{eqn: Dean SDE}, the approximated density for samples of \(\Pi(t)\) yields an approximation to \(P_{\pi}(\cdot, t)\) for \Cref{eqn: Fokker-Planck}. Second, we apply a nonlinear transformation to simplify the support of the density estimation task. In particular, each SDE trajectory \(\Pi(t)\) is supported on a \((d-1)\)-dimensional simplex defined as follows:
\[
\Delta = \{\pi = (\pi_1, \ldots, \pi_{d}) \mid \sum_{i=1}^{d} \pi_i =1, 0 \leq \pi_i \leq 1 \text{ for } i  =1, \ldots, d\}.
\]
To remove the simplex constraint of \(\Delta\), we use a logarithmic transformation to transform the distribution from the \(\pi\) variable to a log-ratio coordinate \(s\). The transformation \(\pi \to s\) ensures that one only needs to apply density estimation over a simple support.

The previous two steps reduce the task to performing density estimation on a lattice model. We use the functional hierarchical tensor over a wavelet basis (FHT-W) ansatz \cite{tang2025wavelet} to perform the density estimation task. The FHT-W model uses a hierarchical tensor model \cite{hackbusch2009new,hackbusch2012tensor} to parameterize density functions. Our numerical experiments in \Cref{sec: numerical experiments} show accurate results for several test cases in 1D and 2D Dean-Kawasaki models. We also obtain accurate results for the challenging setting of \(d = 64\). \edit{Representing the coefficient tensor of the density function without a tensor network model would require a storage complexity that is exponential in \(d\), which is why using a tensor network for dimension reduction is essential.}

\edit{The key idea behind the wavelet coordinate system is to exploit the multi-scale structure of lattice models: after a wavelet transformation \(s \to c\), the distribution \(P_{c}(c, t)\) admits a lower-rank HT representation than the original distribution \(P_{s}(s, t)\) \cite{tang2025wavelet}. The discretized DK model is a lattice model in which adjacent grid cells are strongly correlated due to the diffusion operator, making it well-suited for this wavelet-based approach. While the FHT-W ansatz is used throughout this work, our formulation is compatible with other tree-based tensor network architectures, including functional tensor trains \cite{chen2023combining, ren2023high, dektor2021dynamic, dektor2021rank, soley2021functional} and the HT format applied directly in the original coordinates \cite{hackbusch2009new,tang2024solving}. The primary contribution of this work is to provide a comprehensive computational framework for studying the DK model, rather than the tensor network ansatz itself.}

\subsection{Background in the Dean-Kawasaki model}\label{sec:background}

\paragraph{Interacting particle perspective}
The Dean-Kawasaki derivation in \cite{dean1996langevin} shows that \Cref{eqn: Dean eq} is associated with the following \(N\)-particle stochastic differential equation on \(X_1, \ldots, X_N\). For \(i = 1, \ldots, N\), one writes 
\begin{equation}\label{eqn: Interacting particle langevin}
    dX_i(t) = - \left(\nabla V_1(X_i(t)) +\frac{1}{N}\sum_{j=1}^{N}\nabla V_{2}(X_i(t) - X_j(t))\right) dt + \sqrt{2\beta^{-1}}dB_i(t),
\end{equation}
where \(B_1(t), \ldots, B_N(t)\) are \(N\) independent Brownian motions in \(\R^n\). 
In particular, one sees that \Cref{eqn: Interacting particle langevin} is an interacting particle equation with independent diffusion terms. Generalized versions of the Dean-Kawasaki model have found numerous applications in agent-based modeling \cite{helfmann2021interacting}, continuum modeling \cite{li2019harnessing}, and are used to describe bacterial
dynamics \cite{thompson2011lattice}, reaction-diffusion dynamics \cite{kim2017stochastic}, and nucleation dynamics \cite{lutsko2012dynamical}. 

\paragraph{Gradient flow perspective}
We provide some intuitions on how the discretization is derived through a gradient flow characterization of the Dean-Kawasaki model (see derivation in \cite{dean1996langevin}). The diffusion term and the \(V_1, V_2\) term in the Dean-Kawasaki model correspond to the energy
\begin{equation}\label{eqn: energy}
    \mathcal{E}(\pi) = \frac{1}{\beta}\int \pi(x) \log (\pi(x))\,dx + \int V_1(x) \pi(x)\,dx +  \frac{1}{2}\int \pi(x)V_2(x-y)\pi(y)\,dxdy,
\end{equation}
and the Dean-Kawasaki model can be written as
\begin{equation}\label{eqn: Dean eq 2}
    \partial_t \pi(x,t) = \mathrm{div}\left(\pi\nabla\frac{\delta \mathcal{E}(\pi)}{\delta \pi}\right) + \mathrm{div}\left(\sqrt{\frac{2\pi}{\beta N}}\eta\right),
\end{equation}
where \(\frac{\delta \mathcal{E}(\pi)}{\delta \pi}\) is the Fr\'echet derivative.
One can see that \Cref{eqn: Dean eq 2} indicates that the Dean-Kawasaki model is a stochastic Wasserstein gradient flow. In particular, the \(\bm{\mu}\) term in \Cref{eqn: Dean SDE} is obtained by taking a cell average of the deterministic flux term \(\pi\nabla\frac{\delta \mathcal{E}(\pi)}{\delta \pi}\), and the \(\bm{\sigma}\) term in \Cref{eqn: Dean SDE} is obtained by taking a cell average of the stochastic flux term \(\sqrt{\pi}\eta\). This gradient flow characterization is essential to deriving the discrete Dean-Kawasaki model, and this approach is also taken by \cite{russo2021finite}.

\subsection{Related work}

\paragraph{Density estimation in lattice models}
There is a substantial body of research on density estimation subroutines applicable to lattice models. In the continuous case, notable neural network models include restricted Boltzmann machines (RBMs)~\cite{hinton2010practical, salakhutdinov2010efficient}, energy-based models (EBMs)~\cite{hinton2002training, lecun2006tutorial, gutmann2010noise,marchand2023multiscale}, normalizing flows~\cite{tabak2010density, tabak2013family, rezende2015variational, papamakarios2021normalizing}, variational auto-encoders (VAEs)~\cite{doersch2016tutorial, kingma2019introduction}, diffusion and flow-based models~\cite{sohl2015deep,zhang2018monge,song2019generative,song2020score, song2021maximum, ho2020denoising, albergo2022building, liu2022flow, lipman2022flow, albergo2023stochastic}.

Tensor-network density estimators represent multivariate functions and densities via low-rank tensor formats combined with fixed univariate bases (polynomials, Fourier, etc.); this viewpoint is closely connected to classical decompositions such as tensor train and hierarchical Tucker \cite{oseledets2011tensor,grasedyck2010hierarchical,hackbusch2012tensor}. For high-dimensional Fokker-Planck problems, low-rank tensor methods have also been explored in deterministic PDE solvers (e.g., TT/QTT time-stepping and cross approximation); see \cite{dolgov2012fast,chertkov2021solution} and references therein. For the approximation power of hierarchical tensor representations applied to function classes, we refer to \cite{schneider2014approximation}. 

In this work, we follow a sample-based route: for each time $t$, we approximate the law of the discretized SDE from Monte Carlo trajectories and fit the resulting density by a functional tensor network. A distinctive feature in the DK setting is that the discretized state lives on a probability simplex, which motivates the centered log-ratio transformation used in \Cref{sec: mapping simplex to euclidean}. 
Prior work \cite{chen2023combining, tang2024solving} considered functional tensor networks for Fokker-Planck equations over Euclidean domains; the present work differs in that we solve over a probability simplex and treat higher-dimensional systems ($d = 64$) with external potentials and pairwise interactions. It is also possible to combine tensor networks with neural networks \cite{wang2022tensor} by using a neural network parameterization for the function basis.

\paragraph{The Dean-Kawasaki model}
The Dean-Kawasaki model has long been of research interest in non-equilibrium statistical physics due to its direct modeling of microscopic interacting particle dynamics. The Dean-Kawasaki model is also an example of fluctuating hydrodynamics \cite{hauge1973fluctuating,donev2010accuracy,balboa2012staggered,delong2013temporal,jack2014geometrical}. It is noteworthy that the Dean-Kawasaki model does not admit martingale solutions for generic choices of parameters and initial condition, which was first proven for a special potential-free case in \cite{konarovskyi2019dean} and in \cite{konarovskyi2020dean} for the general case. The well-posedness of the Dean-Kawasaki model does not directly relate to our work, as the spatial discretization is sufficient to provide an implicit regularization to the Dean-Kawasaki model. Further regularization has been proposed to ensure the stability of the Dean-Kawasaki model \cite{cornalba2019regularized,cornalba2023regularised,cornalba2023density}. As our approach is suitable for density estimation over a general simplex, our methodology can be extended to the regularized Dean-Kawasaki model.

\section{Discretization of the Dean-Kawasaki model}\label{sec: dean}
We first describe the approach for spatially discretizing the Dean-Kawasaki model to form a finite-dimensional SDE. For simplicity, we consider an \(n\)-dimensional case with \(n = 1, 2\), and the domain is an \(n\)-dimensional cube \(\Omega = [0, 1]^n\) prescribed with periodic boundary conditions in each physical dimension. The finite-volume treatment generalizes to high-dimensional cases, which we omit for simplicity. Our simple approach is motivated by the finite-volume approaches adopted in \cite{kim2017stochastic, russo2021finite}. Any reasonable finite-volume discretization can be used, and the choice only results in a modification of the associated Fokker-Planck equation in \Cref{eqn: Fokker-Planck}. Therefore, our subsequent tensor network density estimation approach can accommodate other discretization procedures.

We give an outline for this section. \Cref{sec: 1D discretization} provides a discretization of the Dean-Kawasaki model in 1D. In \Cref{sec: 2D discretization}, we present the derivation in 2D, which is more notationally involved and can be skipped on a first reading.

\subsection{1D discretization}\label{sec: 1D discretization}
We give the derivation for the finite-volume discretization to simulate the 1D Dean-Kawasaki model in \Cref{eqn: Dean eq}. The goal is to derive an SDE as shown in \Cref{eqn: Dean SDE}, which would then lead to the Fokker-Planck equation in \Cref{eqn: Fokker-Planck}. 

We take \(h = \frac{1}{m}\) to be the spatial grid mesh size and take \(x_{j} = (j-\frac{1}{2})h \) for \(j = 1, \ldots, m\) to be the cell center points. We take \(x_{j+\frac{1}{2}} = jh\) to be the cell boundary points. The domain is subdivided into grid cells \(C_{j} = [x_{j-\frac{1}{2}}, x_{j+\frac{1}{2}}] = [(j-1)h, jh]\). The finite-volume approach takes a cell-average of \(\pi\) over each \(C_j\). Formally, we write
\[
\Pi(t) = (\Pi_1(t), \ldots, \Pi_m(t)), \quad \Pi_j(t) = \frac{1}{h}\int_{C_j} \pi(x, t) \,dx.
\]

By the discussion in \Cref{sec:background}, we see that the Dean-Kawasaki model is a stochastic gradient flow driven by a deterministic flux term and a stochastic flux term. For \(\mathcal{E}\) as defined in \Cref{eqn: energy}, one has the two flux terms written as follows
\[
F_{d}(x, t) = \pi\nabla\frac{\delta \mathcal{E}(\pi)}{\delta \pi}, \quad  F_{s}(x, t) = \sqrt{\frac{2\pi}{\beta N}}\eta,
\]
where \(\eta\) is a space-time white noise, \(F_d\) is the deterministic flux, and \(F_s\) is the stochastic flux. 

One can then derive the discretized Dean-Kawasaki model by integrating \Cref{eqn: Dean eq 2}. By the divergence theorem, one has 
\[
    d\Pi_j(t) = \frac{1}{h}\left(F_{d, j+\frac{1}{2}} - F_{d, j - \frac{1}{2}}\right)dt + \frac{1}{h}\left(F_{s, j+\frac{1}{2}} - F_{s, j - \frac{1}{2}}\right),
\]
where \(F_{d, j+\frac{1}{2}}, F_{s, j+\frac{1}{2}}\) are the evaluations of \(F_d\) and \(F_s\) at the boundary \(x_{j+\frac{1}{2}}\). In particular, the approximation for \(F_{d, j+\frac{1}{2}}, F_{d, j-\frac{1}{2}}\) forms the \(\bm{\mu}\) term in \Cref{eqn: Dean SDE}, and the approximation for \(F_{s, j+\frac{1}{2}}, F_{s, j-\frac{1}{2}}\) forms the \(\bm{\sigma}\) term in \Cref{eqn: Dean SDE}. 

To further simplify the equation, we use \(F_{V,j+\frac{1}{2}}\) to denote the flux at the boundary \(x_{j+\frac{1}{2}}\) corresponding to the \(V_1, V_2\) terms in \(\mathcal{E}\). For the diffusion term corresponding to the entropy term in \(\mathcal{E}\), we use a simple central difference scheme, leading to an SDE
\begin{equation}\label{eqn: Dean discrete 1D}
\begin{aligned}
    d\Pi_j(t) 
    &=
    \frac{1}{\beta h^2} \left(\Pi_{j+1}(t) + \Pi_{j-1}(t) - 2\Pi_{j}(t)\right)dt\\
    &+
    \frac{1}{h}\left(
    F_{V,j + \frac{1}{2}}
    -
    F_{V,j - \frac{1}{2}}
    \right)dt + \frac{1}{h}\left(F_{s, j+\frac{1}{2}} - F_{s, j - \frac{1}{2}}\right).
\end{aligned}
\end{equation}

First, we describe how one approximates \(F_{V, j+\frac{1}{2}}\). We use the Roe scheme \cite{roe1981approximate}. Define \(F_{V, j}\) to be the flux term through \(x_{j}\) corresponding to the \(V_1, V_2\) terms in \(\mathcal{E}\). The discretization takes the form
\[
F_{V, j} =  \Pi_{j}(t) V_1'(x_{j}) + \frac{1}{m}\sum_{i}V_2'(x_j - x_i)\Pi_{i}(t)\Pi_{j}(t).
\]
Then, for the \(V_1, V_2\) term, the corresponding flux term is defined by
\[
F_{V, j+\frac{1}{2}} = \begin{cases}
    F_{V,j} & A_{V, j+\frac{1}{2}} < 0,\\
    F_{V,j+1} & A_{V, j+\frac{1}{2}} \geq 0,
\end{cases}
\]
where
\[
A_{V, j+\frac{1}{2}} = \frac{F_{V, j+1} - F_{V, j}}{\Pi_{j+1}(t) - \Pi_{j}(t)}.
\]

For the stochastic flux term \(F_{s, j+\frac{1}{2}}\), we take the approach described in \cite{kim2017stochastic}. Let \(B_{j+\frac{1}{2}}(t)\) be an independent Wiener processes for \( j =0, \ldots, m - 1\). The discretization takes the form
\[
F_{s, j+\frac{1}{2}} = \sqrt{\frac{2\tilde{\Pi}_{j+\frac{1}{2}}(t)}{h\beta N}}dB_{j+ \frac{1}{2}}(t),
\]
where \(\tilde{\Pi}_{j+\frac{1}{2}}\) is the spatial average of \(\pi(x, t)\) around \(x = x_{j+\frac{1}{2}}\). In particular, \(\tilde{\Pi}_{j+\frac{1}{2}}\)
is a modified arithmetic mean defined by
\[
\tilde{\Pi}_{j+\frac{1}{2}}(t) = \frac{\Pi_j(t) + \Pi_{j+1}(t)}{2}H\left(h\Pi_j(t)\right)H\left(h\Pi_{j+1}(t)\right),
\]
where \(H\) is a Heaviside function defined as:
\[
H(x) = \begin{cases}
    0 & x \leq 0,\\
    x & 0 \leq x \leq 1,\\
    1 & x \geq 1.
\end{cases}
\]
The choice of \(F_{s, j + \frac{1}{2}}\) follows \cite{kim2017stochastic} and has the advantage of improved numerical stability. When \(V_1 = V_2 = 0\), it has been shown in \cite{kim2017stochastic} that the resulting continuous-time SDE in \Cref{eqn: Dean discrete 1D} preserves the positivity of \(\Pi_j(t)\) for any \(j, t\).

The aforementioned definition of \(F_{V, j+\frac{1}{2}}\) and \(F_{s, j+\frac{1}{2}}\) fully characterizes the SDE in \Cref{eqn: Dean discrete 1D}. Overall, one sees that the finite-volume discretization enables the SPDE in \Cref{eqn: Dean eq} to be written as a tractable SDE in \Cref{eqn: Dean discrete 1D}. To simulate the SDE, we perform an additional temporal discretization. Let \(T\) be the terminal time. We take \(\delta t\) to be a fixed time step, and we assume that \(\delta t\) is an integer fraction of \(T\) with \(K = T/\delta t\). We define \(t_{k} = k\delta t\) for \(k = 0, \ldots, K\). We take \(\Pi^{k}_{j}\) to model \(\Pi_{j}(t_k)\). 

For the initial condition, we set \(\Pi^{0}_{j} = \frac{1}{h}\int_{C_j}\pi(x, 0) \, dx\); One can approximate the integral through numerical integration or a simple point evaluation. For the flux term, we take \(F^{k}_{V, j + \frac{1}{2}}, F^{k}_{s, j + \frac{1}{2}}\) to approximate the value of \(F_{V, j + \frac{1}{2}}\) and \(F_{s, j + \frac{1}{2}}\) at time \(t = t_k\). Subsequently, with the identification \(\Pi^{k}_{0} = \Pi^{k}_{m}\) and \(\Pi^{k}_{m+1} = \Pi^{k}_{1}\), we use the following scheme for \(j = 1, \ldots, m\):
\begin{equation}\label{eqn: Euler-Maruyama 1D}
\begin{aligned}
    \Pi_{j}^{k+1}  
    =
    &\Pi_{j}^{k}
    +
    \frac{\delta t}{\beta h^2} \left(\Pi^{k+1}_{j+1} + \Pi^{k+1}_{j-1} - 2\Pi^{k+1}_{j}\right)
    \\
    +
    &\frac{\delta t}{h}\left(
    F^{k}_{V,j + \frac{1}{2}}
    -
    F^{k}_{V,j - \frac{1}{2}}
    \right)
    +
    \frac{\sqrt{\delta t}}{h}\left(
    F^{k}_{s,j + \frac{1}{2}}
    -
    F^{k}_{s,j - \frac{1}{2}}
    \right),
\end{aligned}
\end{equation}
where we use an implicit scheme for the diffusion term \(\Delta \pi\) to ensure stability under larger \(d t\).

Thus, it remains to fully specify \(F^{k}_{V, j + \frac{1}{2}}\) and \(F^{k}_{s, j + \frac{1}{2}}\). We first specify \(F^{k}_{V, j + \frac{1}{2}}\). 
We use \(F_{V, j}^{k}\) to represent \(F_{V, j}\) at time \(t_k\). We write
\[
 F^{k}_{V, j} =  \Pi_{j}^{k} V_1'(x_{j}) + h\Pi_{j}^{k}\sum_{i}V_2'(x_j - x_i)\Pi_{i}^{k}.
\]
Then, we define \(F^{k}_{V, j + \frac{1}{2}}\) by the following:
\[
F^{k}_{V,j+\frac{1}{2}} = \begin{cases}
    F^{k}_{V,j} & A^{k }_{V,j+\frac{1}{2}} < 0,\\
    F^{k}_{V,j+1} & A^{k}_{V,j+\frac{1}{2}} \geq 0,
\end{cases} \quad A^{k}_{V,j+\frac{1}{2}} = \frac{F_{V,j+1}^{k} - F_{V,j}^{k}}{\Pi_{j+1}^{k} - \Pi_{j}^{k}}.
\]

We then specify \(F^{k}_{s, j + \frac{1}{2}}\).
Let \(W^{k}_{j+ \frac{1}{2}} \sim N(0, 1)\) be i.i.d. standard normal random variables for \(k = 0, \ldots, K, j = 0, \ldots, m-1\). 
We write
\[
F_{s, j+ \frac{1}{2}}^{k} = \sqrt{\frac{2\tilde{\Pi}^{k}_{j+\frac{1}{2}}}{h\beta N}}W_{j+ \frac{1}{2}}^{k}, \quad
\tilde{\Pi}^{k}_{j+\frac{1}{2}} =  \frac{\Pi_j^k + \Pi_{j+1}^k}{2}H\left(hN\Pi_j^k\right)H\left(hN\Pi_{j+1}^k\right).
\]

Thus, simulation with respect to the SDE in \Cref{eqn: Dean discrete 1D} allows one to simulate the discretized Dean-Kawasaki model. For any \(t = t_k\), the simulation results of \((\Pi^{k}_{j})_{j = 1}^{m}\) form an approximation to the SPDE solution \(\pi(\cdot, t_k)\). In this case the dimension is \(d = m\). For practical purposes, we perform a clamping on the Gaussian variable \(W_{j+ \frac{1}{2}}^{k}\). We replace \(W_{j+ \frac{1}{2}}^{k}\) with \(\max(-5, \min(W_{j+ \frac{1}{2}}^{k}, 5))\) to prevent blow-up from rare large normal draws. In practice, a more careful treatment with noise clamping is left for future work.
For the density estimation, we shall overload the notation \(\pi\) and work with the associated cell-mass vector
\[
\pi(t) := h\,\Pi(t),
\]
which satisfies $\pi_j(t)\ge 0$ and (for a mass-conserving discretization) $\sum_{j=1}^m \pi_j(t)=1$.
This is the simplex variable used in \Cref{sec: mapping simplex to euclidean}. Likewise, in the 2D case, we use \(\pi(t) := h^2\,\Pi(t)\).

\subsection{2D discretization}\label{sec: 2D discretization}
A similar procedure applies to 2D, with modest modifications from the 1D case. For simplicity, we only define the spatial discretization of the SDE. The temporal discretization follows from a simple modification of \Cref{eqn: Euler-Maruyama 1D}. We take \(h = \frac{1}{m}\) and take \(x_{j} = (j-\frac{1}{2})h \) for \(j = 1, \ldots, m\). We take \(x_{j+\frac{1}{2}} = jh\). Let \(x_{(i,j)}=  (x_i, x_j)\). The domain is subdivided into grid cells \(C_{(i,j)} = [x_{i-\frac{1}{2}}, x_{i+\frac{1}{2}}] \times [x_{j-\frac{1}{2}}, x_{j+\frac{1}{2}}]   = [(i-1)h, ih] \times [(j-1)h, jh]\). The finite-volume approach takes a cell average of \(\pi\) over each \(C_{i,j}\). Formally, we write
\[
\Pi(t) = (\Pi_{(1,1)}(t), \ldots, \Pi_{(m,m)}(t)), \quad \Pi_{(i,j)}(t) = \frac{1}{h^2}\int_{C_{(i,j)}} \pi(x, t) \,dx,
\]
and \(\Pi(t)\) is a \(d\)-dimensional discretization of \(\pi\) with \(d = m^2\).

In particular, the gradient flow perspective in \Cref{sec: 1D discretization} carries through. One writes
\[
(F_{d,1}(x, t), F_{d,2}(x, t)) = \pi\nabla\frac{\delta \mathcal{E}(\pi)}{\delta \pi}, \quad  (F_{s,1}(x, t), F_{s,2}(x, t)) = \sqrt{\frac{2\pi}{\beta N}}\eta,
\]
and the divergence theorem leads to
\[
\begin{aligned}
    d\Pi_{(i,j)}(t) = &
    \frac{1}{h}\left(F_{d,1, (i+\frac{1}{2}, j)} - F_{d,1, (i-\frac{1}{2}, j)}\right)dt +\frac{1}{h}\left(F_{s,1, (i+\frac{1}{2}, j)} - F_{s,1, (i-\frac{1}{2}, j)}\right)\\
    +&\frac{1}{h}\left(F_{d,2, (i, j+\frac{1}{2})} - F_{d,2, (i, j-\frac{1}{2})}\right)dt +\frac{1}{h}\left(F_{s,2, (i, j+\frac{1}{2})} - F_{s,2, (i, j-\frac{1}{2})}\right),
\end{aligned}
\]
where \(F_{d,1, (i+\frac{1}{2}, j)}, F_{s,1, (i+\frac{1}{2}, j)}\) are respectively the average flux of \(F_{d,1}, F_{s,1}\) into the right boundary of \(C_{(i,j)}\), and \(F_{d,2, (i, j+\frac{1}{2})}, F_{s,2, (i, j+\frac{1}{2})}\) are respectively the average flux of \(F_{d,2}, F_{s,2}\) into the upper boundary of \(C_{(i,j)}\). To formally define \(F_{d,1,(i+\frac{1}{2}, j)}\), for example, one writes
\[
F_{d,1,(i+\frac{1}{2}, j)} = \frac{1}{h}\int_{x_{j-\frac{1}{2}}}^{x_{j+\frac{1}{2}}}F_{d,1}( x_{i+\frac{1}{2}}, z)\,dz,
\]
and the rest of the formal definition follows similarly.

Similar to the 1D case, we separate the deterministic flux term by using a central difference scheme to approximate the diffusion term. We use \(F_{V,1, (i+\frac{1}{2}, j)}, F_{V,2, (i, j+\frac{1}{2})}\) to represent the part of \(F_{d,1, (i+\frac{1}{2}, j)}, F_{d,2, (i, j+\frac{1}{2})}\) corresponding to the \(V_1, V_2\) term. Then, in this case, the SDE reads
\begin{equation}\label{eqn: Dean discrete 2D}
\begin{aligned}
    &d\Pi_{(i,j)}(t)\\ = 
    &\frac{1}{\beta h^2} \left(\Pi_{(i+1,j)}(t)  + \Pi_{(i-1,j)}(t)   + \Pi_{(i,j+1)}(t)   + \Pi_{(i,j-1)}(t)   - 4\Pi_{(i,j)}(t)  \right)dt
    \\
    +&\frac{1}{h}\left(F_{V,1, (i+\frac{1}{2}, j)} - F_{V,1, (i-\frac{1}{2}, j)}\right)dt + \frac{1}{h}\left(F_{s,1, (i+\frac{1}{2}, j)} - F_{s,1, (i-\frac{1}{2}, j)}\right)\\
    +&\frac{1}{h}\left(F_{V,2, (i, j+\frac{1}{2})} - F_{V,2, (i, j-\frac{1}{2})}\right)dt +\frac{1}{h}\left(F_{s,2, (i, j+\frac{1}{2})} - F_{s,2, (i, j-\frac{1}{2})}\right).
\end{aligned}
\end{equation}

First, we describe how one approximates \(F_{V, l, (i+\frac{1}{2},j+\frac{1}{2})}\). We again use a Roe scheme. Define \(F_{V, l, (i,j)}\) to be the flux term through \(x_{(i,j)}\) corresponding to the \(V_1, V_2\) term in \(\mathcal{E}\). The discretization takes the following form for \(l = 1, 2\):
\[
\begin{aligned}
    F_{V, l, (i,j)} = &\Pi_{(i,j)}(t) \partial_{x_l}V_1(x_{(i,j)})+ \frac{1}{d}\sum_{i',j'}\partial_{x_l}V_2(x_{(i,j)} - x_{(i',j')})\Pi_{(i',j')}(t)\Pi_{(i,j)}(t).
\end{aligned}
\]
Then, for \(l = 1, 2\), the numerical flux is defined by
\begin{equation*}
    \begin{aligned}
        &F_{V, 1, (i + \frac{1}{2}, j)} = \begin{cases}
        F_{V, 1, (i,j)} & A_{V, 1, (i+\frac{1}{2},j)} < 0,\\
        F_{V, 1, (i+1,j)} & A_{V, 1, (i+\frac{1}{2},j)} \geq 0,
    \end{cases}
    \quad 
    A_{V,1, (i+\frac{1}{2},j)} = \frac{F_{V,1,(i+1,j)} - F_{V,1,(i,j)}}{\Pi_{(i+1,j)}(t) - \Pi_{(i,j)}(t)}.
    \\
        &F_{V, 2, (i, j+\frac{1}{2})} = \begin{cases}
        F_{V, 2, (i,j)} & A_{V, 2, (i,j+\frac{1}{2})} < 0,\\
        F_{V, 2, (i,j+1)} & A_{V, 2, (i, j+\frac{1}{2})} \geq 0,
    \end{cases}
    \quad 
    A_{V,2,(i,j+\frac{1}{2})} = \frac{F_{V,2,(i,j+1)} - F_{V,2,(i,j)}}{\Pi_{(i,j+1)}(t) - \Pi_{(i,j)}(t)}.
    \end{aligned}    
\end{equation*}

For the stochastic flux terms \(F_{s, 1, (i+\frac{1}{2}, j)}, F_{s, 2, (i, j+\frac{1}{2})}\), we likewise use the scheme in \cite{kim2017stochastic}. Let \(B_{l, (i+\frac{1}{2},j+\frac{1}{2})}(t)\) be independent Wiener processes for \( i,j =0, \ldots, m-1\) and \(l = 1,2\). Then, the discretization takes the form
\[
F_{s,l,(i + \frac{1}{2}, j)} = \sqrt{\frac{2\tilde{\Pi}_{l, ( i+\frac{1}{2}, j+\frac{1}{2})}(t)}{h^2\beta N}}dB_{l,( i+\frac{1}{2}, j+\frac{1}{2})}(t),
\]
where \(\tilde{\Pi}_{1, ( i+\frac{1}{2}, j+\frac{1}{2})}(t)\) and \(\tilde{\Pi}_{2, ( i+\frac{1}{2}, j+\frac{1}{2})}(t)\) are both spatial approximations of \(\pi\) around
\((x_{i+\frac{1}{2}}, x_{j+\frac{1}{2}})\).

Subsequently, one can perform similar time discretizations as in \Cref{sec: 1D discretization}, which allows one to simulate the Dean-Kawasaki model. We skip the implementation details for brevity. The detailed procedures for simulation can be found in the GitHub repository accompanying this manuscript. Simulating the SDE in \Cref{eqn: Dean discrete 2D} allows one to simulate the discretized Dean-Kawasaki model. For any \(t = t_k\), the simulation results of \((\Pi^{k}_{(i,j)})_{i,j = 1}^{m}\) form an approximation to the SPDE solution \(\pi(\cdot, t_k)\). For the 2D case, we also use clamping to ensure stability.

\section{Main formulation}
This section details the approach for performing density estimation over a \((d-1)\)-dimensional simplex \(\Delta\) with a functional hierarchical tensor. Our goal is to solve the Fokker-Planck equation in \Cref{eqn: Fokker-Planck} with density estimation. Thus, the input to this subroutine is a collection of samples on \(\Delta\) at time \(t\), and the output is a normalized density \(P_{\pi}(\pi, t)\) in a functional tensor network ansatz.

Explicitly, by simulating the SDE in \Cref{eqn: Dean discrete 1D} and \Cref{eqn: Dean discrete 2D}, we obtain a collection of samples. For the case of 1D, we obtain a collection of \(B\) samples \(\{\pi^{(b)}(t) = (h\Pi^{(b)}_{j}(t))_{j = 1, \ldots, m}\}_{ b = 1, \ldots, B}\). Likewise, in the case of 2D, we obtain a collection of \(B\) samples \(\{\pi^{(b)}(t) = (h^2\Pi^{(b)}_{(i,j)}(t))_{i,j = 1, \ldots, m}\}_{\, b = 1, \ldots, B}\). 
The samples at time \(t\) form an empirical distribution approximation to the Fokker-Planck equation solution \(P_{\pi}(\pi, t) \approx \frac{1}{B}\sum_{b = 1}^{B}\delta(\pi - \pi^{(b)}(t))\).
In other words, the discretized Dean-Kawasaki model introduced in \Cref{sec: dean} describes a specific probability distribution over \(\Delta\), and performing the density estimation over the simulated samples can give a particle-based approximation to the associated Fokker-Planck equation in \Cref{eqn: Fokker-Planck}.

We outline this section by expanding on how we address a few challenges in the aforementioned density estimation approach. First, conventional tensor network density estimation only applies to Euclidean spaces, and the simplex constraint in \Cref{eqn: simplex} for the Fokker-Planck equation poses challenges for tensor networks. To address the issue, we use an entry-wise logarithmic transformation followed by an entry-wise offset to center the distribution. Overall, this coordinate transformation maps \(\Delta\) to \(\R^{d}\) and effectively eliminates the simplex constraints. This coordinate transformation is covered in the first part of \Cref{sec: mapping simplex to euclidean}.

Following the logarithmic coordinate transformation, one obtains a lattice model, and a challenge is that the lattice model exhibits strong site-wise coupling for practically interesting models. To efficiently address the strong site-wise coupling, we perform a wavelet transformation that separates the model into different length scales. The coordinate transformation to the wavelet coefficient is covered in the second part of \Cref{sec: mapping simplex to euclidean}. After the transformation to the wavelet coefficient space, we use the functional hierarchical tensor over a wavelet basis (FHT-W) ansatz for density estimation. The FHT-W model is covered in \Cref{sec: ansatz}.

Lastly, the transformation into wavelet coefficients creates challenges for accurate observable estimation. Essentially, once the density function in the form of an FHT-W ansatz is obtained, we also require an FHT-W representation of the observable. We propose using interpolation to compress the observable function in the FHT-W ansatz. The interpolation procedure is covered in \Cref{sec: interpolation}, which is technical and can be skipped for a first reading.

\paragraph{Notation}
For notational compactness, we introduce several shorthand notations for subsequent derivations. For \(n \in \mathbb{N}\), let \([n] := \{1,\ldots, n\}\). For an index set \(S \subset [d]\), we let \(x_{S}\) stand for the subvector with entries from the index set $S$.

\subsection{Mapping \texorpdfstring{$\Delta$}{simplex} to Euclidean space}\label{sec: mapping simplex to euclidean}
We describe an invertible transformation of \(\Delta\) to \(\R^{d-1}\) through a coordinate transformation from \(\pi\) to \(c\). This procedure consists of two coordinate transformations \(\pi \to s\) followed by \(s \to c\), which we cover in sequence.

\paragraph{Transform to the log space}
The \((d-1)\)-dimensional simplex \(\Delta\) is defined as follows:
\begin{equation}\label{eqn: simplex}
        \Delta = \{\pi = (\pi_1, \ldots, \pi_{d}) \mid \sum_{i=1}^{d} \pi_i = 1, 0 < \pi_i \leq 1 \text{ for } i  =1, \ldots, d\},
\end{equation}
and we see that \(\Delta \subset \R^{d}\) is a convex polytope. It is less suitable to directly perform density estimation of \(\Delta\) in the \(\pi\) coordinate. Instead, we propose applying a nonlinear transformation to \(\Delta\) as a subset of \(\R^{d}\). We define the map \(s \colon \Delta \to \R^{d}\) by the following formula:
\begin{equation}\label{eqn: transform}
    s(\pi_1, \ldots, \pi_{d}) = (\log(\pi_1) - Y, \ldots, \log(\pi_{d}) - Y), \quad Y = \frac{1}{d}\sum_{j = 1}^{d}\log(\pi_{j}).
\end{equation}

The transformation above is a centered log-ratio map: it takes $\pi$ in the interior of the simplex to a vector $s\in\R^d$ satisfying $\sum_{i=1}^d s_i=0$.
In particular, the uniform state $\pi = (1/d, \ldots, 1/d)$ maps to $s= (0, \ldots, 0)$.

The inverse map is the softmax:
\[
\pi_i(s)=\frac{\exp(s_i)}{\sum_{j=1}^d \exp(s_j)}.
\]
Thus, the interior of $\Delta$ is diffeomorphic to the hyperplane $\{s\in\R^d:\sum_i s_i=0\}$.

For concreteness, consider $d=3$. The simplex is $\Delta=\{(\pi_1,\pi_2,\pi_3)\in\R_{\geq 0}^3:\pi_1+\pi_2+\pi_3=1\}$, an open triangle.
The centered log-ratio map sends $\pi$ to $s\in\R^3$ with $s_1+s_2+s_3=0$ via
\[
s_i=\log(\pi_i)-\tfrac13(\log\pi_1+\log\pi_2+\log\pi_3).
\]
The uniform point $(1/3,1/3,1/3)$ maps to $s=(0,0,0)$.
When $\pi_1\to 0$ (with $\pi_2,\pi_3$ held fixed), $s_1\to -\infty$; hence the boundary of $\Delta$ is mapped to infinity in the $s$-plane.
In computations, we restrict to a truncated simplex $\Delta_\varepsilon=\{\pi\in\Delta:\min_i\pi_i\ge\varepsilon\}$, whose image is bounded. This is consistent with the stable-regime diagnostics described above.

We explain how \Cref{eqn: transform} is designed. The formula for \(s\) consists of two steps: (a) perform entrywise logarithm to get \(\log(\pi)\), (b) subtract the mean to get \(\log(\pi) - Y\). We write \(s(\pi) = s  = (s_1, \ldots, s_d)\). This nonlinear transformation transforms \(\Delta\) from the \(\pi\) coordinate to the \(s\) coordinate. One can see that \(s(\pi)\) is the centered log-ratio coordinate of \(\pi\). The map \(\pi \to (s_1, \ldots, s_{d-1})\) is a parameterization of \(\Delta\) into an open subset of \(\R^{d-1}\). 
In practice, the numerical simulation of the SDE in \Cref{eqn: Dean SDE} is unstable if any entry of \(\Pi(t)\) is too close to zero. \edit{For example, the noise term in \Cref{eqn: Dean discrete 1D} is proportional to $\sqrt{\tilde{\Pi}_{j+1/2}}$, whereas the diffusion term is proportional to $(\Pi_{j+1} + \Pi_{j-1} - 2\Pi_j)$.} In a region where density is low, the noise effect is larger than the diffusion effect, which makes it challenging for the simulation to remain stable. Therefore, in the fluctuating hydrodynamics literature, it is standard to assume that the diffusion term is strong (see Assumption FD2, FD4 of \cite{cornalba2023dean}). In other words, for any sample point \(\Pi(t) = (\pi_1, \ldots, \pi_d)\) generated from the procedure in \Cref{sec: dean}, the diffusion is sufficiently strong that \(\pi_k\) is bounded away from zero with high probability. \edit{Therefore, under the strong-diffusion (stable-regime) assumption, the induced \(s\)-samples concentrate in a bounded region of the hyperplane \(\sum_i s_i = 0\). Additionally, for robustness, if occasional near-boundary events occur, we apply a standard positivity-preserving clamping step at the SDE level (see \Cref{sec: dean}).
Under this stable-regime assumption, the induced $s$-samples concentrate in a bounded region of the hyperplane $\sum_i s_i=0$, and it becomes meaningful to approximate the density on a bounded computational domain after an additional normalization step (see below).}

We see that the $\pi$ to $s$ transformation already allows a standard density estimation procedure (i.e., one operating on a Euclidean domain without simplex constraints) to be applied. Let \(P_{\pi}(\pi, t)\) denote the target distribution over \(\Delta\) under the \(\pi\) coordinate. With the \(s\) transformation, we can instead consider \(P_{s}(s, t)\), which is the target distribution under the \(s\) coordinate. Formally, \(P_{s}(s, t)\) is the pushforward of \(P_{\pi}(\pi, t)\) under \(s(\cdot)\). 
Moreover, one is given samples from \(P_{s}(s, t)\) by transforming the samples of \(P_{\pi}(\pi, t)\) through \(s(\cdot)\). One can directly apply the FHT ansatz in \cite{tang2024solving} to perform density estimation from samples of \(P_{s}(s, t)\).

\paragraph{Transform to the wavelet space}
We perform a further coordinate transformation through a wavelet transform.
In our discretization for the Dean-Kawasaki model, the distribution \(P_{\pi}(\pi, t)\) is a lattice model in 1D and 2D. With the pushforward into the log space, \(P_{s}(s, t)\) is also a lattice model. 
\edit{In our experiments, we see that the distribution $P_s(s, t)$ exhibits strong correlations between nearby components $s_j$ and $s_{j+1}$. The wavelet transformation $s \to c$ separates the distribution into different length scales and reduces the numerical rank of the resulting distribution $P_c(c, t)$. This rank reduction has been demonstrated in \cite{tang2025wavelet} for several lattice models with strong nearest-neighbor correlations.} Moreover, the constraint that entries of \(s\) sum to zero can be elegantly addressed by the wavelet transformation.

We first go through the iterative wavelet coarsening procedure in the 1D case. For simplicity, we assume that \(d = 2^{L}\) so that \(\pi = (\pi_1, \ldots, \pi_{2^L})\).
In the first step of the coarse-graining process, the wavelet filter transforms the \(s\) variable into two \(2^{L-1}\)-dimensional variables \(y_{L - 1} = (y_{j, L-1})_{j \in [2^{L-1}]}\) and \(c_{L - 1} = (c_{j, L-1})_{j \in [2^{L-1}]}\) defined as
\[
y_{j, L-1} = \frac{s_{2j - 1} + s_{2j}}{\sqrt{2}}, \quad c_{j, L-1} = \frac{s_{2j - 1} - s_{2j}}{\sqrt{2}}.
\]
The variable \(y_{L-1}\) is the scaling coefficient of \(s\) at length scale \(2^{-(L-1)}\), and \(c_{L-1}\) is the detail coefficient of \(s\) at length scale \(2^{-(L-1)}\).

Repeating the same approach, for \(l = L-1, \ldots, 1\), one applies the wavelet transform to \(y_{l}\) to obtain \(2^{l-1}\)-dimensional variables \(y_{l-1}\) and \(c_{l-1}\) given by the following equation
\[
y_{j, l-1} = \frac{y_{2j - 1, l} + y_{2j, l}}{\sqrt{2}}, \quad c_{j, l-1} = \frac{y_{2j - 1, l} - y_{2j, l}}{\sqrt{2}}.
\]
Finally, we note that \(y_1\) is \(2\)-dimensional and the procedure ends at \(y_1 \to (y_0, c_0)\). As \(\sum_{j = 1}^{2^L}s_j = 0\) due to \Cref{eqn: transform}, we have \(y_{0} = 0\). Therefore, we have fully described the procedure to transform \(s = (s_1, \ldots, s_d)\) into \((c_{k, l})_{k, l}\). In particular, the transformation from \(\pi\) to \(s = (s_1, \ldots, s_d)\) to \((c_{k, l})_{l = 0, \ldots, L - 1, \, k = 1, \ldots, 2^l}\) is an invertible embedding from \(\Delta\) to a subset of \(\R^{d-1}\).

By the same procedure, one can perform the iterative wavelet transformation on the \(s\) variable in the 2D case. One can perform a sequence of iterative 2D wavelet transforms to the \(s\) variable \(s = (s_{(i,j)})_{i,j = 1}^{m}\) with \(m = 2^{L/2}\). By the procedure in \cite{tang2025wavelet}, one likewise transforms the \(s\) variable to \((c_{k, l})_{l = 0, \ldots, L - 1, \, k = 1, \ldots, 2^l}\). The map from \(\pi\) to \(c\) is also invertible in this case.

Thus, for time \(t \in [0, T]\), one is given a collection of the \(B\) samples \(\{\Pi^{(b)}(t) = (\Pi^{(b)}_{j}(t))_{j = 1, \ldots, m}\}_{b = 1, \ldots, B}\) in 1D or \(\{\Pi^{(b)}(t) = (\Pi^{(b)}_{(i,j)}(t))_{i,j = 1, \ldots, m}\}_{ b = 1, \ldots, B}\) in 2D. Applying the map \(\pi \to s\) followed by the 1D/2D wavelet transformation yields samples in the \(c\) space. We denote the resulting wavelet-coefficient samples by \(\{C^{(b)}(t) = (C^{(b)}_{k, l}(t))_{l = 0, \ldots, L - 1, \, k = 1, \ldots, 2^l}\}_{ b = 1, \ldots, B}\). The samples \(\{C^{(b)}(t)\}_{ b = 1, \ldots, B}\) are the inputs to the subsequent density estimation subroutine. Moreover, the samples form an empirical distribution approximation to \(P_{c}(c, t)\). We subsequently use a tensor network ansatz to perform density estimation over the samples.

\subsection{Tensor network architecture of \texorpdfstring{$P_{c}(c, t)$}{Pc(c, t)}}\label{sec: ansatz}

\begin{figure}[h]
    \centering
    \includegraphics[width=
    \linewidth]{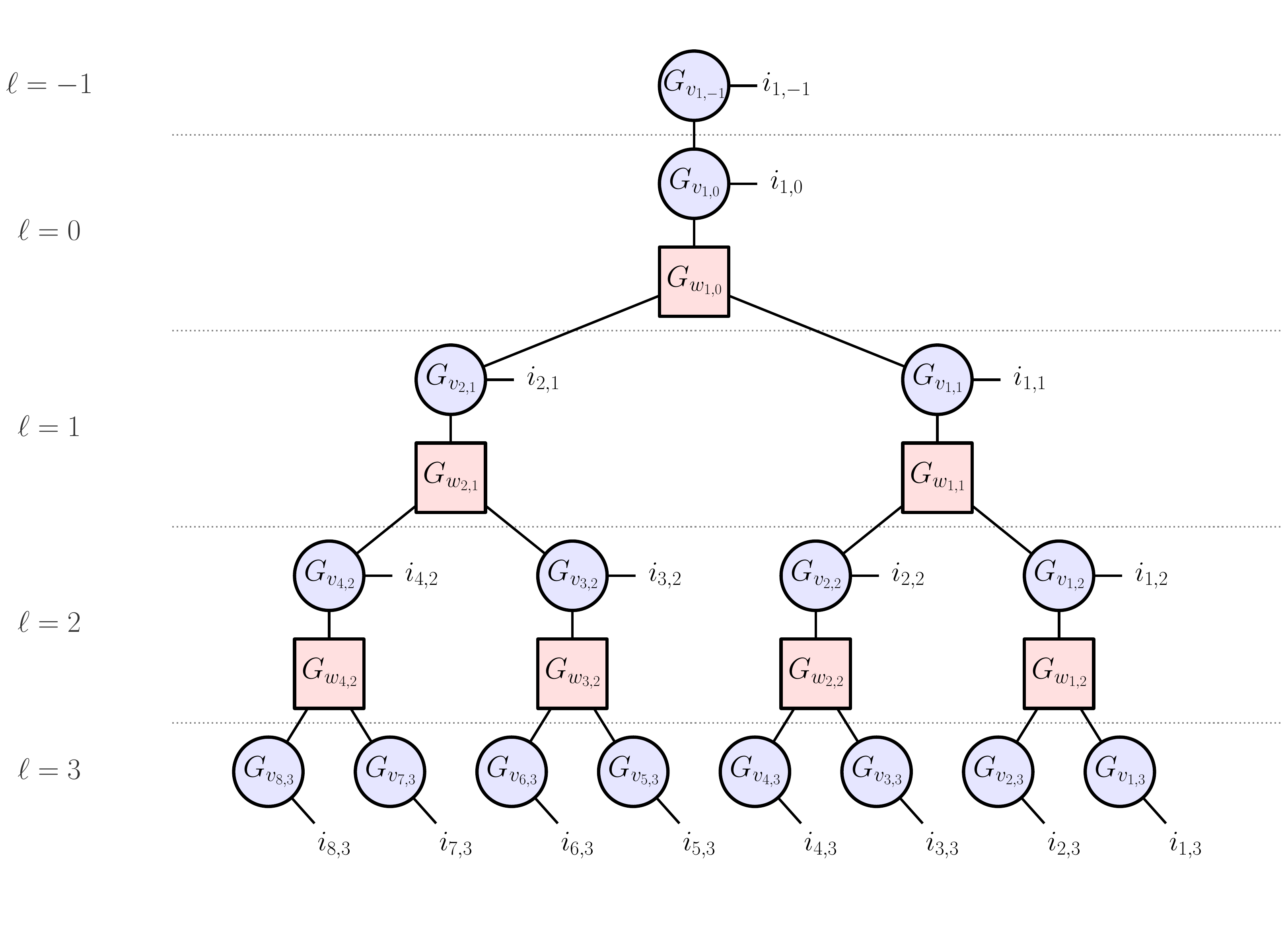}
    \caption{Tensor diagram of the FHT-W ansatz for \(L = 4\)
    (\(d = 2^L = 16\) wavelet coordinates).
    Circles (blue): tensor components \(G_{v_{k,l}}\) at the external
    nodes \(v_{k,l} \in V_{\mathrm{ext}}\); each carries one physical
    bond \(i_{k,l}\) of dimension \(q+1\), indexing a Legendre
    polynomial basis for the wavelet coordinate \(c_{k,l}\).
    Squares (pink): tensor components \(G_{w_{k,l}}\) at the internal
    nodes \(w_{k,l} \in V \setminus V_{\mathrm{ext}}\); these have no
    physical bonds.
    Edges between nodes are contracted bonds of dimension at most \(r\).
    Dashed lines separate wavelet resolution levels from the coarsest
    (\(l = -1\)) to the finest (\(l  = 3\)).
    The diagram represents the coefficient tensor \(D\) in
    \Cref{eq:ttn-contraction_sec_TTN_with_internal};
    in the experiments \(q \in \{15, 25\}\) and \(r = 20\)
    (see \Cref{tab:parameters}).}
    \label{fig:wavelet_FHT_L_3}
\end{figure}

After the transformation from the \(\pi\) coordinate into the
wavelet coordinate \((c_{k, l})_{k, l}\), we use a functional
hierarchical tensor over a wavelet basis (FHT-W) ansatz to model
the distribution \(P_{c}(c, t)\). The index \(l\) represents the
resolution scale, where a larger \(l\) corresponds to finer-scale
information, and \(k = 1, \ldots, 2^{l}\) is the spatial location
index at level \(l\). \edit{The FHT-W model is a tree tensor network
(hierarchical Tucker format) whose binary tree mirrors the dyadic
structure of the wavelet transform: at each level \(l\), an
internal node \(w_{k,l}\) connects the wavelet coordinate
\(c_{k,l}\) to two finer-level coordinates \(c_{2k-1,l+1}\) and
\(c_{2k,l+1}\), as illustrated in \Cref{fig:wavelet_FHT_L_3}.}

\begin{remark}[Relationship to the hierarchical Tucker format]
The HT decomposition was introduced by Hackbusch and K\"uhn \cite{hackbusch2009new} for general tensor spaces, including function spaces represented with respect to a basis system; see also \cite{hackbusch2012tensor,grasedyck2010hierarchical}. The FHT-W ansatz used in this work is an instance of the HT format in which: (a) a wavelet coordinate system is used to improve the low-rank approximability of lattice models with nearest-neighbor correlations, and (b) a collection of non-leaf nodes contains a physical bond. We use the acronym FHT-W following the tensor network terminology in \cite{tang2024solving,tang2025wavelet} to distinguish the specific wavelet-adapted variant. We emphasize that the underlying decomposition is the HT format of \cite{hackbusch2009new}.
\end{remark}

\edit{The internal nodes \(w_{k,l}\) reduce the tensor dimensions at each branch point, keeping the overall contraction cost linear in the number of coordinates. A full description of the FHT-W architecture and density estimation procedure can be found in \cite{tang2025wavelet}.}

For readers unfamiliar with tree-based functional tensor networks, we briefly summarize the background needed to understand FHT-W. Details can be found in \cite{tang2025wavelet}.

\paragraph{Notations for tree graphs}
A tree graph \(T = (V, E)\) is a connected undirected graph without cycles. An undirected edge \(\{v,v'\} \in E\) is written interchangeably as $(v, v')$ or $(v', v)$. For any \(v \in V\), let \(\mathcal{N}(v)\) denote the neighbors of \(v\). Moreover, let \(\mathcal{E}(v)\) denote the set of edges incident to \(v\). For an edge \(e = (k, w)\), removing the edge \(e\) in \(E\) results in two connected components \(I_{1}, I_{2}\) with \(I_1 \cup I_2 = V\). 
For any \((w, k) \in E\), we use \(k \to w\) to denote the unique connected component among \(I_1, I_2\) which contains \(k\) as a node.

\begin{figure}
    \centering
    \subfloat[\centering]{\includegraphics[width=0.45\linewidth]{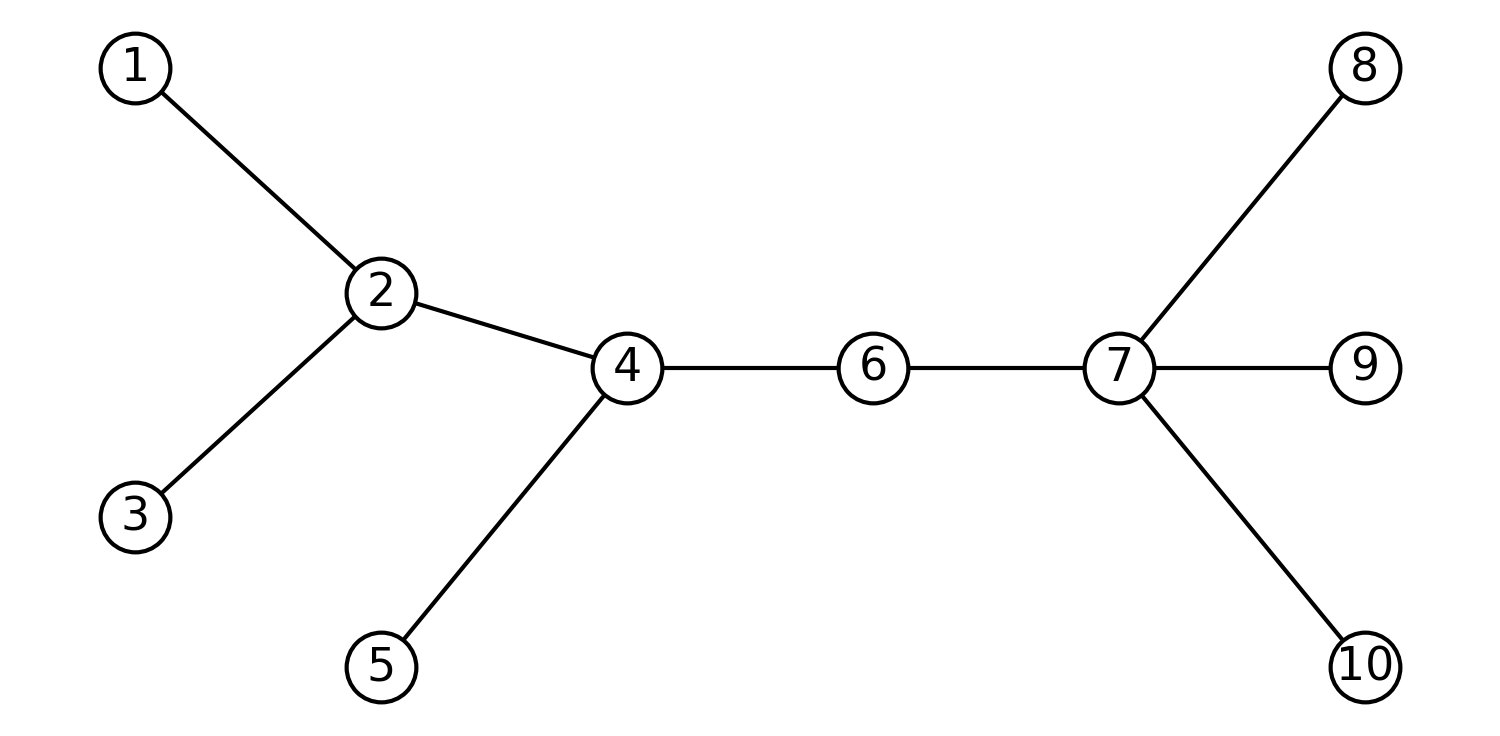}}%
    \hspace{12pt}
    \subfloat[\centering]{\includegraphics[width=0.45\linewidth]{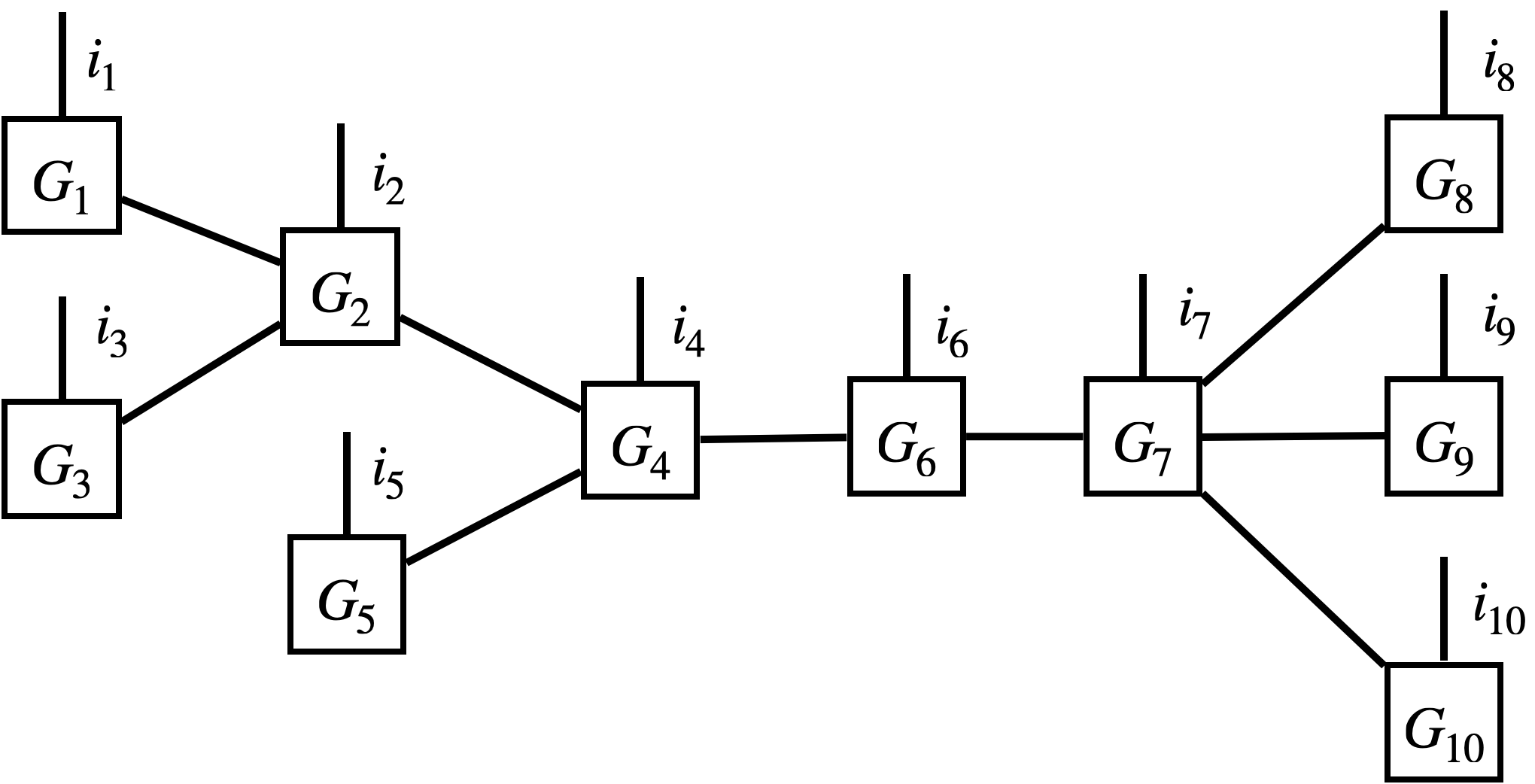}}%
    \caption{(A) A tree structure \(T = (V, E)\) with \(V = V_{\mathrm{ext}}= \{1, \ldots, 10\}\). (B) Tensor Diagram representation of a tree tensor network over \(T\).}
    \label{fig:tree+ttn_SEC_TTN}
\end{figure}

\paragraph{Tree-based FTN}
We introduce the tree-based functional tensor network ansatz for density estimation. Let \(T = (V, E)\) be a tree graph and let \(V_{\mathrm{ext}}\) be the subset of vertices in \(T\) whose associated tensor component admits an external bond. To adhere to the setting of \(d\)-dimensional functions, let \(|V_{\mathrm{ext}}| = d\). The total vertex size \(\tilde{d} = |V|\) is equal to \(d\) only when the tensor network has no internal nodes. We define the tree-based functional tensor network in \Cref{def: tree-based FTN}. The tensor network is illustrated in \Cref{fig:tree+ttn_SEC_TTN}.
\begin{definition}\label{def: tree-based FTN}
    (Tree-based functional tensor network)
    Suppose one has a tree structure $T = (V, E)$ and the corresponding ranks $\{r_e: e \in E\}$. The symbol \(V_{\mathrm{ext}}\) is the subset of \(V\) that contains external bonds. For simplicity, we apply a labeling of nodes so that \(V = [\tilde{d}]\) and \(V_{\mathrm{ext}} = [d]\).
    
    The \emph{tensor component} at \(k \in V\) is denoted by \(G_{k}\). Let \(\mathrm{deg}(k)\) stand for the degree of \(k\) in \(T\).
    When \(k \in V_{\mathrm{ext}}\), \(G_{k}\) is defined as a \((\mathrm{deg}(k) + 1)\)-tensor of the following shape:
    \begin{equation*}\label{eqn: core size constraint case external}
        G_{k} : [n_{k}] \times \prod_{e \in \mathcal{E}(k)}  [r_e] \rightarrow \R.
    \end{equation*}
    When \(k \not \in V_{\mathrm{ext}}\), \(G_{k}\) is defined as a \(\mathrm{deg}(k)\)-tensor of the following shape:
    \begin{equation*}\label{eqn: core size constraint case internal}
        G_{k} : \prod_{e \in \mathcal{E}(k)}  [r_e] \rightarrow \R.
    \end{equation*}
    
    A \(d\)-tensor \(D\) is said to be a \emph{tree tensor network} defined over the tree \(T\) and tensor components $\{G_k\}_{k=1}^{\tilde{d}}$ if
    \begin{equation}
        \label{eq:ttn-contraction_sec_TTN_with_internal}
        D(i_1, \ldots, i_{d}) = \sum_{\alpha_E} \prod_{k 
        =1}^{d} G_k\left(i_k, \alpha_{\mathcal{E}(k)}\right)\prod_{k = d+1}^{\tilde{d}}G_k\left( \alpha_{\mathcal{E}(k)}\right).
    \end{equation}

    We use $n_k$ to denote the basis size at node $k$; this should not be confused with $n \in \{1,2\}$, the spatial dimension of the PDE domain introduced in \Cref{sec:intro}.
    For \(j = 1, \ldots, d\), we let \(\{\psi_{i,j}\}_{i = 1}^{n_j}\) denote a collection of orthonormal univariate basis functions for the \(j\)-th coordinate. A tree-based FTN is the \(d\)-dimensional function defined by \(D\), and it admits the following equation:
    \begin{equation}
        \label{eq: tree-based FTN forward map}
        p(x_1, \ldots, x_{d}) = \sum_{ \alpha_E} \prod_{k 
        =1}^{d} \left( \sum_{i_{k} = 1}^{n_k}G_k\left(i_k, \alpha_{\mathcal{E}(k)}\right)\psi_{i_k, k}(x_k)\right)\prod_{k = d+1}^{\tilde{d}}G_k\left( \alpha_{\mathcal{E}(k)}\right).
    \end{equation}
\end{definition}

In particular, in the experiments of \Cref{sec: numerical experiments}, we set $n_k = q + 1$ for all physical nodes (i.e., nodes in $V_{\mathrm{ext}}$), where $q$ is the maximal polynomial degree.

\paragraph{Tree structure for FHT-W}

We first specify the vertex of \(T\) in the case of the FHT-W model used for the discretized Dean-Kawasaki model SDE. We remark that the data structure of the \(c = (c_{k, l})_{k,l}\) variable is the same for both the 1D and 2D cases. For the vertex set, we define the external nodes to be \(V_{\mathrm{ext}} = \{v_{k, l}\}_{l \in [L-1] \cup \{0\}, k \in [2^{\max(l, 0)}]}\). The total vertex set is \(V = V_{\mathrm{ext}} \cup \{w_{k, l}\}_{l \in [L-2]\cup \{0\}, k \in [2^l]}\). The vertex \(v_{k, l}\) corresponds to the variable \(c_{k, l}\), and the variables \(w_{k, l}\) is an internal node inserted at each level \(l = 0, \ldots, L-2\).

The edge set of \(T\) is \(E = E_1 \cup E_2\), where \(E_1 = \{(w_{k, l}, v_{k, l})\}_{w_{k, l} \in V}\) and \(E_2 = \{(w_{k, l}, v_{2k-1, l+1}), (w_{k, l}, v_{2k, l+1})\}_{w_{k, l} \in V}\). The definition of \(E\) shows that each variable \(c_{2k-1, l+1}\) and \(c_{2k, l+1}\) is connected to the variable \(c_{k, l}\) through an internal node \(w_{k, l}\). This construction for \(T\) ensures that \(c_{k, l}\) is placed close to \(c_{k', l'}\) when \((k, l)\) is close to \((k', l')\). The description fully characterizes the tree \(T = (V, E)\). One sees that the description aligns with the illustration in \Cref{fig:wavelet_FHT_L_3}.

\subsection{Interpolation on tree-based FTN}
\label{sec: interpolation}
The goal of this subsection is to prepare a target function \(M(c)\) in the FHT-W ansatz. Assuming this is done, the observable estimation task of estimating \(\mathbb{E}_{C \sim P_{c}(c, t)}\left[M(C)\right]\) can be performed by an efficient tensor diagram contraction between \(P_{c}(c, t)\) and \(M(c)\). The notation for the algorithm is dense, which is necessary to accommodate different topologies that might arise in the tree-based FTN. For understanding, it might be instructive to refer to \cite{tang2024solvingb}, where the interpolation procedure is introduced for a notationally simpler case of \(T\) being a binary tree with no internal nodes. 

\paragraph{Relation to black-box tensor reconstruction and cross interpolation}
The interpolation procedure used here is closely related to black-box reconstruction of low-rank tensors in hierarchical Tucker/tree tensor formats and to cross-interpolation type methods.
In particular, the idea of recovering cores from point evaluations by selecting informative slices appears in the black-box HT literature \cite{ballani2013black} and in more recent maxvol-based HT reconstruction methods \cite{ryzhakov2024black}. In our work, recovery is still done by point evaluation. However, as we are reconstructing a functional tensor network, our interpolation is done by performing informative measurements on the HT rather than selecting informative slices.

For future reference, we introduce the interpolation scheme for tree-based functional tensor networks. In this procedure, the input is the target function \(M\), which has an analytic expression, and the output is an FHT or FHT-W ansatz that approximates \(M\). As both the FHT and the FHT-W ansatz might be used for the Dean-Kawasaki model, we introduce a general interpolation subroutine under the general tree structure. 
\paragraph{Example of observable function}
We provide an example of an observable function that is compressed in practice.
Suppose one has used the FHT-W model to approximate \(P_{c}(c, t)\). An important observable estimation task is based on the observable function \[M(\pi) = -\sum_{j = 1}^{d}\pi_j\log(\pi_j),\] which is the entropy of \(\pi\). Performing observable estimation over \(M(\pi)\) is important as it measures the disorder in the system. The entropy term also appears naturally in the variational form of the Fokker-Planck equation \cite{jordan1998variational} and the Dean-Kawasaki model \cite{jack2014geometrical}. Let \(M(c)\) be the pushforward of \(M(\pi)\) with \(M(c) = M(\pi(c))\). To be rigorous, let \(s^{-1}(\cdot)\) denote the transformation from \(s\) to the \(\pi\) variable, and let \(\mathrm{IDWT}(\cdot)\) be the inverse wavelet transformation, and then one can write \(M(c) = M(s^{-1}(\mathrm{IDWT}(c)))\).
One can see that the function \(M(c)\) has a complicated functional form and thus can only be prepared in an FHT-W through interpolation.

The derivation in \cite{tang2024solvingb} shows that interpolation and density estimation in the FHT ansatz can be done by solving structured independent linear equations. The FHT ansatz is a special case of the tree-based FTN, and we show that the interpolation of a function \(M(x)\) in a tree-based FTN can be done similarly by solving linear equations.

\paragraph{Linear equation for \(G_k\)}

We assume that \(M(x)\) is a tree-based FTN defined over a tree graph \(T = (V, E)\). Our goal is to solve for each \(G_k\) for \(k \in V\) by a tractable linear equation.
We first derive the linear equation for each \(G_{k}\). Let \(e_1, \ldots, e_{\mathrm{deg}(k)}\) be all edges incident to \(k\), and let \(e_{j} = (v_j, k)\) for \(j = 1, \ldots, \mathrm{deg}(k)\). For each \(v_{j}\), we construct \(D_{v_{j} \to k}\) by contracting all tensor components \(G_{w}\) for \(w \in v_{j} \to k\). When \(k \in V_{\mathrm{ext}} = [d]\), \Cref{def: tree-based FTN} implies the following linear equation:
\begin{equation}\label{eqn: core-determining equation}
    D(i_1, \ldots, i_d) = \sum_{a_{\mathcal{E}(k)}}G_{k}(i_k, a_{\mathcal{E}(k)}) \prod_{j = 1}^{\mathrm{deg(k)}}D_{v_j \to k}(\alpha_{e_j},i_{v_j \to k \cap [d]}),
\end{equation}
and the equation for \(k = d + 1 ,\ldots, \tilde{d}\) can be obtained from \Cref{eqn: core-determining equation} by simply omitting the \(i_{k}\) index in the right-hand side of the equation. 

\paragraph{Sketched equation for \(G_k\)}
We see that \Cref{eqn: core-determining equation} is exponential-sized and is thus not practical. Overall, one solves for \(G_k\) by solving a sketched linear equation.
We first consider the case for \(k \in V_{\mathrm{ext}} = [d]\).
For \(j = 1, \ldots, \mathrm{deg}(k)\), we let \(\tilde{r}_{e_{j}} > r_{e_{j}}\) be some integer, and we introduce sketch tensors \(S_{v_j \to k} \colon \prod_{l \in v_j \to k \cap [d] }[n_l] \times [\tilde{r}_{e_j}]\). By contracting \Cref{eqn: core-determining equation} with \(\bigotimes_{j = 1}^{\mathrm{deg}(k)}S_{v_j \to k}\), one obtains the following sketched linear equation for \(G_k\):
\begin{equation}\label{eqn: TTN sketched linear equation}
    B_{k}(i_k, \beta_{\mathcal{E}(k)}) = \sum_{\alpha_{\mathcal{E}(k)}}\prod_{j = 1}^{\mathrm{deg}(k)}A_{v_j \to k}(\beta_{e_j}, \alpha_{e_j})G_k(i_k, \alpha_{\mathcal{E}(k)}),
\end{equation}
where \(B_{k}\) is the contraction of \(D\) with \(\bigotimes_{j = 1}^{\mathrm{deg}(k)}S_{v_j \to k}\), and each \(A_{v_j \to k}\) is the contraction of \(D_{v_j \to k}\) with \(S_{v_j \to k}\). The equation for \(k=d+1,\ldots, \tilde{d}\) can likewise be obtained from \Cref{eqn: TTN sketched linear equation} by simply omitting the \(i_{k}\) index on both sides.

\paragraph{The choice of sketch tensor}
It suffices to specify the sketch tensors \(S_{v \to v'}\). We assume for concreteness that \(M \colon [-1 ,1]^d \to \R\). 
For each edge \((v, k)\), we select \((v \to k \cap [d])\)-indexed points \(\{w^{\beta}_{v \to k \cap [d]} = (w^{\beta}_{j})_{j \in v \to k \cap [d]} \}_{\beta \in [\tilde{r}_{(v, k)}]} \subset [-1, 1]^{|v \to k| \cap [d]}\).
Corresponding to the selected points, we construct sketch tensor \(S_{v \to k} \colon \prod_{j \in v \to k \cap [d]} [n_{j}] \times [\tilde{r}_{(v, k)}] \to \R\) to be defined by 
\begin{equation}\label{eqn: def of S_V_k}
    S_{v \to k}(:,\beta) = \bigotimes_{j \in v \to k \cap [d]}\Vec{\Psi}_{j}\left(w_{j}^{\beta}\right).
\end{equation}

Moreover, for each node \(k = 1 ,\ldots, d\), we select points \(\{w_k^{\mu_k}\}_{\mu_k \in [\tilde{n}_k]} \subset [-1, 1]\) for \(\tilde{n}_k \geq n_k\). Corresponding to the selected points, we construct the sketch tensor \(S_{k} \colon [\tilde{n}_k] \times [n_k] \to \R\) to be defined by 
\begin{equation}
    S_{k}(\mu_k,i_k) = \psi_{i_k; k}(w_k^{\mu_k}).
\end{equation}

\paragraph{Obtain \(B_k\) by \(M\) queries}
The FTN evaluation of \(M\) defined in terms of \(D\) is written as follows:
\begin{equation}
\label{eqn: FTN forward map}
    M(x_{1}, \ldots, x_{d}) = \sum_{i_{k} \in [n_k], k \in [d]} D_{i_1,\ldots, i_{d}} \psi_{i_1, 1}(x_1)\cdots \psi_{i_{d}, d}(x_{d}) = \left<D, \, \bigotimes_{j=1}^{d} \Vec{\Psi}_{j}(x_j) \right>.
\end{equation}
To obtain \(B_k\), one contracts \(D\) with \(\otimes_{j = 1}^{\mathrm{deg}(k)}S_{v_j \to k}\).
Write \(\beta_j = \beta_{(w_j, k)}\) for \(j = 1, \ldots, \mathrm{deg}(k)\). 
When \(k = 1 ,\ldots , d\), combining \Cref{eqn: def of S_V_k} and \Cref{eqn: FTN forward map} implies
\begin{equation}
\label{eqn: TTNS B_k interpolation formula}
    \sum_{i_k}S_{k}(\mu_k, i_k)B_k(i_k, \beta_{\mathcal{E}(k)}) = M(w^{\beta_{1}}_{v_1 \to k \cap [d]}, \ldots, w^{\beta_{\mathrm{deg}(k)}}_{v_{\mathrm{deg}(k)} \to k \cap [d]}, w^{\mu_k}_{k}),
\end{equation}
and so one can obtain \(B_k\) as the unique least-squares solution to \Cref{eqn: TTNS B_k interpolation formula}. In practice, we directly use the pseudoinverse formula
\begin{equation}
\label{eqn: Tree FTN B_k interpolation formula internal node case external}
    B_k(i_k, \beta_{\mathcal{E}(k)}) = \sum_{\mu_k}S_{k}^{\dagger}(i_k, \mu_k)M(w^{\beta_{1}}_{(v_1 \to k) \cap [d]}, \ldots, w^{\beta_{\mathrm{deg}(k)}}_{(v_{\mathrm{deg}(k)} \to k) \cap [d]}, w^{\mu_k}_{k}).
\end{equation}

When \(k = d + 1 ,\ldots, \tilde{d}\), directly using \Cref{eqn: def of S_V_k} and \Cref{eqn: FTN forward map} implies
\begin{equation}
\label{eqn: Tree FTN B_k interpolation formula internal node case internal}
    B_k(\beta_{\mathcal{E}(k)}) = M(w^{\beta_{1}}_{(v_1 \to k) \cap [d]}, \ldots, w^{\beta_{\mathrm{deg}(k)}}_{(v_{\mathrm{deg}(k)} \to k) \cap [d]}).
\end{equation}

\paragraph{Obtain \(A_k\) by \(M\) queries + SVD}

To obtain \(A_{(k, v)}\), we take the approach in \cite{tang2025wavelet} whereby we first obtain \(Z_{(k, v)}\) from contracting \(C\) with \(S_{k \to v} \otimes S_{v \to k}\). The contraction can be done by querying \(M\), as \Cref{eqn: FTN forward map} and \Cref{eqn: def of S_V_k} imply the equation
\begin{equation}\label{eqn: Z_k_v interpolation formula}
    Z_{(k, v)}(\beta, \gamma) = M(w^\beta_{k \to v}, w^\gamma_{v \to k}).
\end{equation}

One then obtains \(A_{k \to v}, A_{v \to k}\) through SVD of \(Z_{(k, v)}\), resulting in \(Z_{(k, v)}(\beta, \gamma) = \sum_{\mu}U(\beta, \mu)W(\mu, \gamma)\). The choice of \(A_{k \to v} = U\) and \(A_{v \to k}= W\) forms a unique and consistent choice of gauge between the pair \((D_{v \to k}, D_{k \to v})\). As obtaining \(U, W\) requires merging the middle factor of the SVD to either the left or right factor, one can use the root information to specify \(U, W\) fully. In particular, if \(v\) is the parent to \(k\), then one takes \(U\) to be the left orthogonal factor in the SVD. If \(k\) is the parent of \(v\), then one takes \(W\) to be the right orthogonal factor in the SVD.
Finally, one obtains \(A_{k}\) by tensor products of all \(A_{v \to k}\) for \(v \in \mathcal{N}(k)\).

\paragraph{Algorithmic summary}
We summarize the tree-based FTN interpolation algorithm in \Cref{alg:tree-based FTN interpolation internal}.

\begin{algorithm}[h]
\caption{Tree-based FTN interpolation.}
\label{alg:tree-based FTN interpolation internal}
\begin{algorithmic}[1]
\REQUIRE Access to \(M \colon [-1, 1]^d \to \R\).
\REQUIRE Tree \(T = (V, E)\) with \(V = [\tilde{d}]\) and \(V_{\mathrm{ext}} = [d]\). 
\REQUIRE Chosen variable-dependent function basis \(\{\Vec{\Psi}_{j}\}_{j \in [d]}\).
\REQUIRE Collection of sketch tensors \(\{S_{v \to k}, S_{k \to v}\}\) and target internal ranks $\{r_{(k, v)}\}$ for each edge \((k, v) \in E\).
\FOR{each edge \((k, v)\) in \(T\)}
    \STATE Obtain \(Z_{(k,v)}\) by \Cref{eqn: Z_k_v interpolation formula}.
    \STATE Obtain \(A_{k \to v}\) as the left factor of the best rank \(r_{(k, v)}\) factorization of \(Z_{(k,v)}\)
    \STATE Obtain \(A_{v \to k}\) as the right factor of the best rank \(r_{(k, v)}\) factorization of \(Z_{(k,v)}\).
\ENDFOR
\FOR{each node \(k\) in \(T\)}
\IF{node \(k\) contains a physical bond}
    \STATE Obtain \(B_{k}\) by \Cref{eqn: Tree FTN B_k interpolation formula internal node case external}.
\ELSE
    \STATE Obtain \(B_{k}\) by \Cref{eqn: Tree FTN B_k interpolation formula internal node case internal}.
\ENDIF
    \STATE With \(\mathcal{E}(k) = \{(v_j, k)\}_{j = 1}^{\mathrm{deg}(k)}\), collect each \(A_{v_j \to k}\).
    \STATE Obtain \(G_{k}\) by solving the overdetermined linear system \((\bigotimes_{j = 1}^{\mathrm{deg}(k)} A_{v_j \to k})G_{k} = B_{k}\).
\ENDFOR
\end{algorithmic}
\end{algorithm}

\section{Numerical experiments}\label{sec: numerical experiments}

We present the numerical performance of the proposed method for the 1D and 2D Dean-Kawasaki models. \Cref{sec: protocol} gives the experiment details. \Cref{sec: 1D numerical experiment} covers the results for 1D. \Cref{sec: 2D numerical experiment} covers the results for 2D.

\subsection{Reproducibility details and experimental protocol}\label{sec: protocol}
We collect all simulation and approximation parameters in \Cref{tab:parameters} and summarize the numerical choices needed to reproduce the figures:
\begin{itemize}
\item \textbf{Discretization:} grid size $m$ (1D) or $m\times m$ (2D), with $d=m$ or $d=m^2$.
\item \textbf{Simulation:} time step $d t$, horizon $T$, number of trajectories $B$, and clamping rule (see \Cref{sec: dean}).
We report the clamping rate and the minimum mass statistic $\pi_{\min}(t)$ to certify stable-regime simulation.
\item \textbf{Transformations:} centered log-ratio $\pi\mapsto s$; Haar DWT $s\mapsto y$; and affine normalization $y\mapsto c\in[-0.9,0.9]^d\subset[-1,1]^d$ (see \Cref{sec: mapping simplex to euclidean}).
\item \textbf{Basis and ranks:} orthonormal Legendre basis of degree $q$ (so $n_k=q+1$ for each node \(k\)); hierarchical ranks $\{r_e\}$ (we report whether ranks are uniform or level-dependent).
\end{itemize}

\begin{table}[h]
\centering
\caption{Experimental parameters. All experiments use $d = 64$ grid cells ($L = 6$ wavelet levels), inverse temperature $\beta = 0.05$, effective particle count $N = 1000$, terminal time $T = 1$, and bond dimension $r = 20$. For observable interpolation: degree $q = 6$, rank $r = 5$. A sensitivity analysis of $r$ and $q$ is provided in the supplementary material.}
\label{tab:parameters}
\begin{tabular}{@{}lcc@{}}
\toprule
Setting & 1D & 2D \\
\midrule
Spatial dimension $n$ & 1 & 2 \\
Grid cells per dimension $m$ & 64 & 8 \\
Sample paths $B$ & 6000 & 12\,000 \\
$d t$ (no potential) & 0.005 & 0.001 \\
$d t$ (with potential) & 0.0002 & 0.0003 \\
Polynomial degree $q$ (density) & 25 & 15 \\
\bottomrule
\end{tabular}
\end{table}

\paragraph{Error metrics}
For observable interpolation, we take a collection of test points \((c^{(b)})_{b = 1}^{N_{\mathrm{test}}}\), and we report the mean relative prediction error (MRPE)
\begin{equation}
    \mathrm{MRPE} := \frac{\sum_{b=1}^{N_{\mathrm{test}}} \left| \widehat{M}(c^{(b)})-M(c^{(b)})\right|}{\sum_{b=1}^{N_{\mathrm{test}}} \left| M(c^{(b)})\right|}.
\end{equation}

\subsection{1D Dean-Kawasaki models}\label{sec: 1D numerical experiment}

\paragraph{Dean-Kawasaki without potential function}
We first consider a simple model by taking \(V_1 = V_2 = 0\) in \Cref{eqn: Dean eq}. In this case, there is no external field potential or pairwise interaction. This simple case has been of research interest and has been studied in \cite{cornalba2023density}. The model reads
\[
\partial_t \pi = \frac{1}{\beta}\Delta \pi + \mathrm{div}(\sqrt{\frac{2\pi}{\beta N}}\eta).
\]
We apply the sketching-based density estimation algorithm of \cite{tang2025wavelet} to the wavelet-transformed samples, using the parameters listed in \Cref{tab:parameters}.

\begin{figure}
  \centering
  \includegraphics[width = \linewidth]{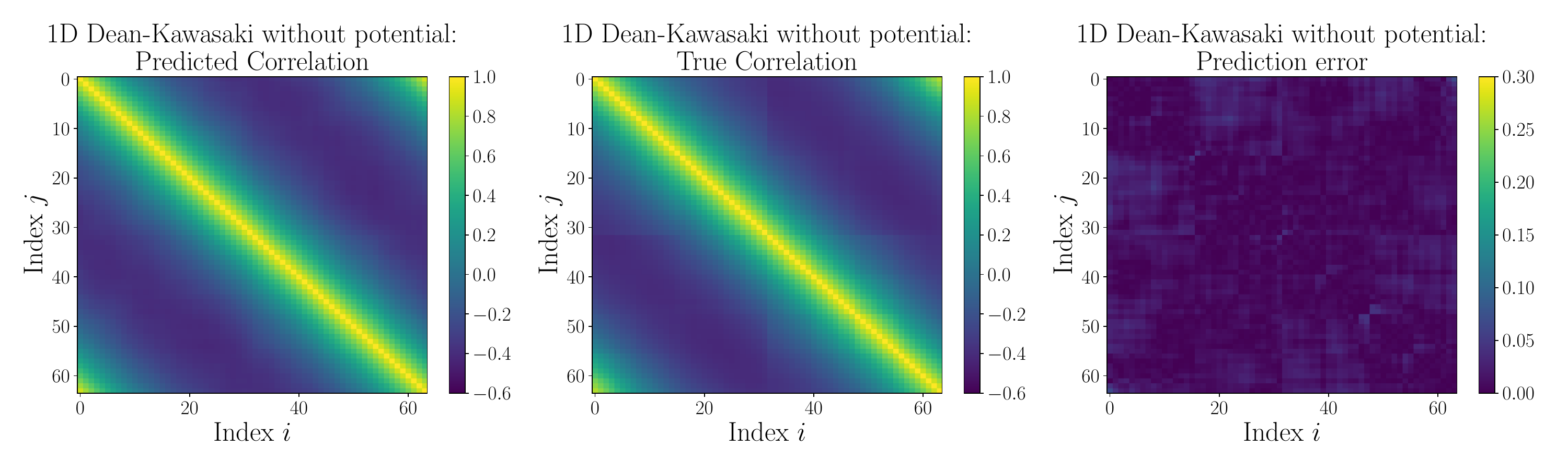}
  \caption{1D Dean-Kawasaki model with \(V_1 = V_2 = 0\). Plot of the correlation matrix predicted by the FHT-W ansatz. The mean prediction error is \(0.011\) and the max prediction error is \(0.060\).}
  \label{Fig: 1D covariance}
  \centering
    \includegraphics[width = \linewidth]{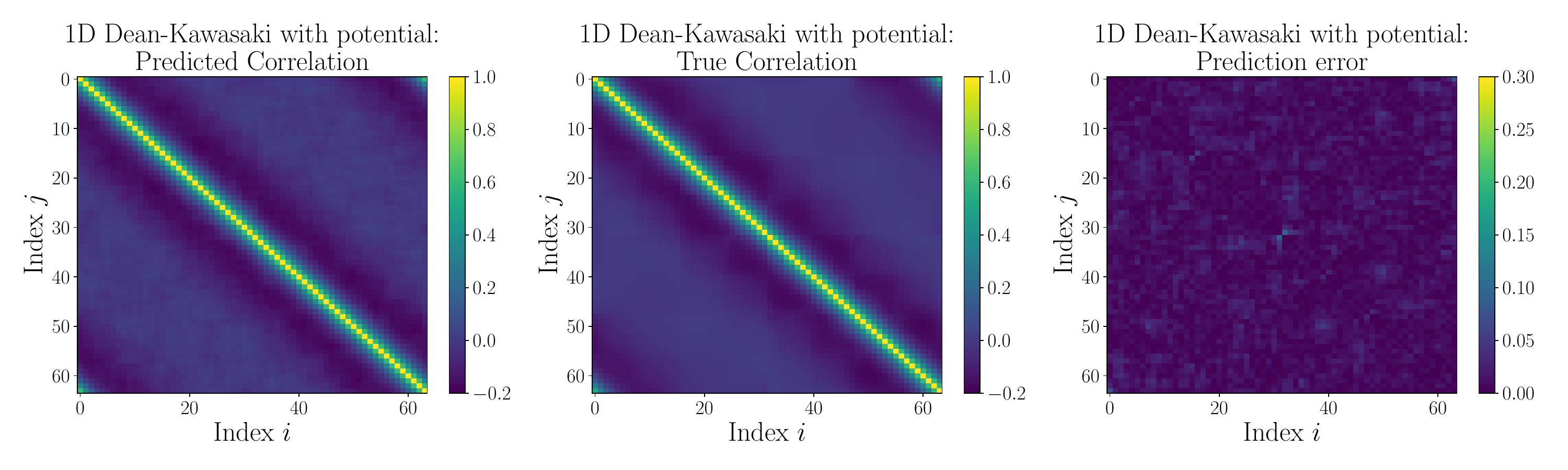}
  \caption{1D Dean-Kawasaki model with an external potential \(V_1\) and a pairwise potential \(V_2\) modeling particle repulsion. Plot of the correlation matrix predicted by the FHT-W ansatz. The mean prediction error is \(0.011\) and the max prediction error is \(0.109\).}
  \label{Fig: 1D covariance Dean potential}
\end{figure}

\begin{figure}[h]
    \centering
    \includegraphics[width= 0.45
    \linewidth]{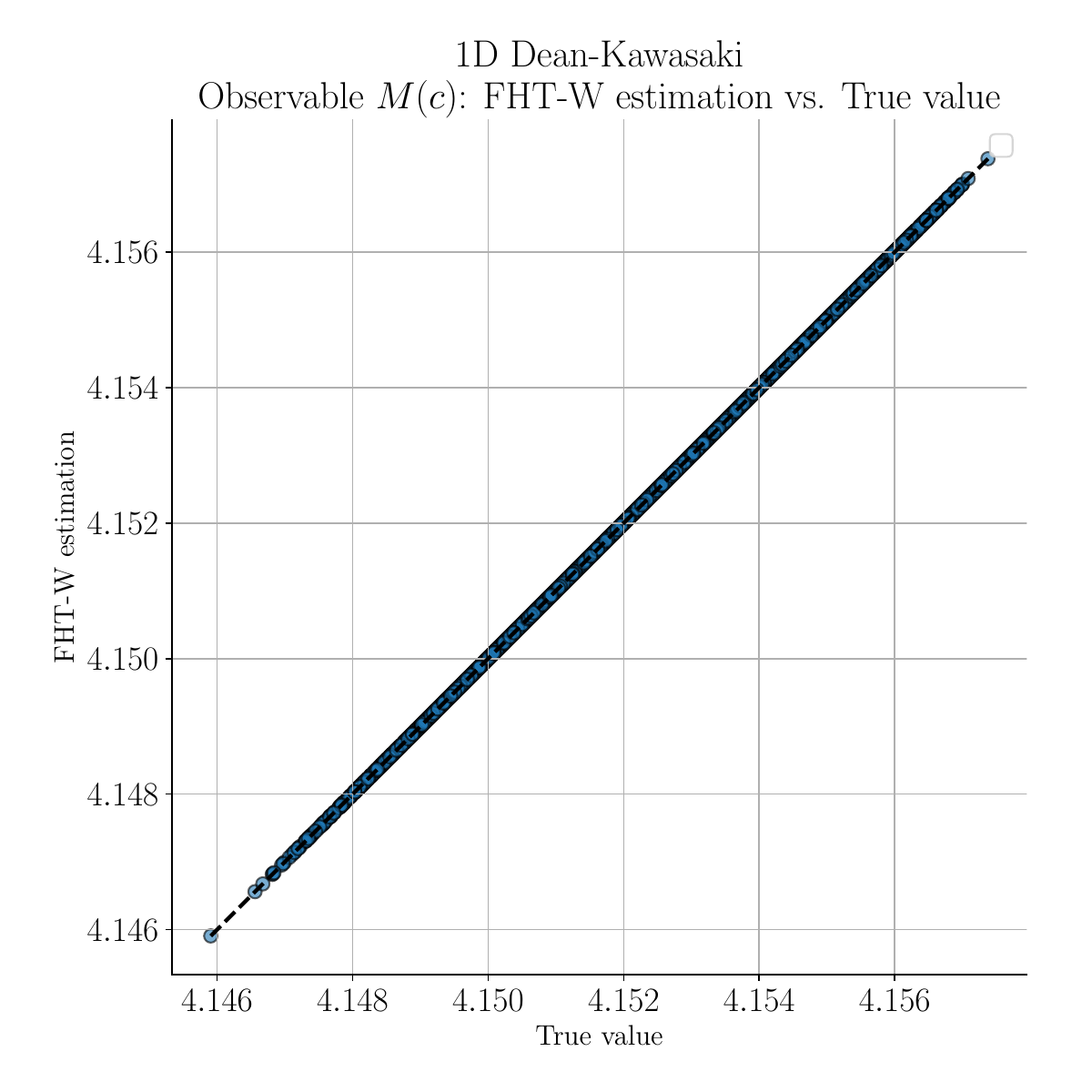}
    \includegraphics[width= 0.45
    \linewidth]{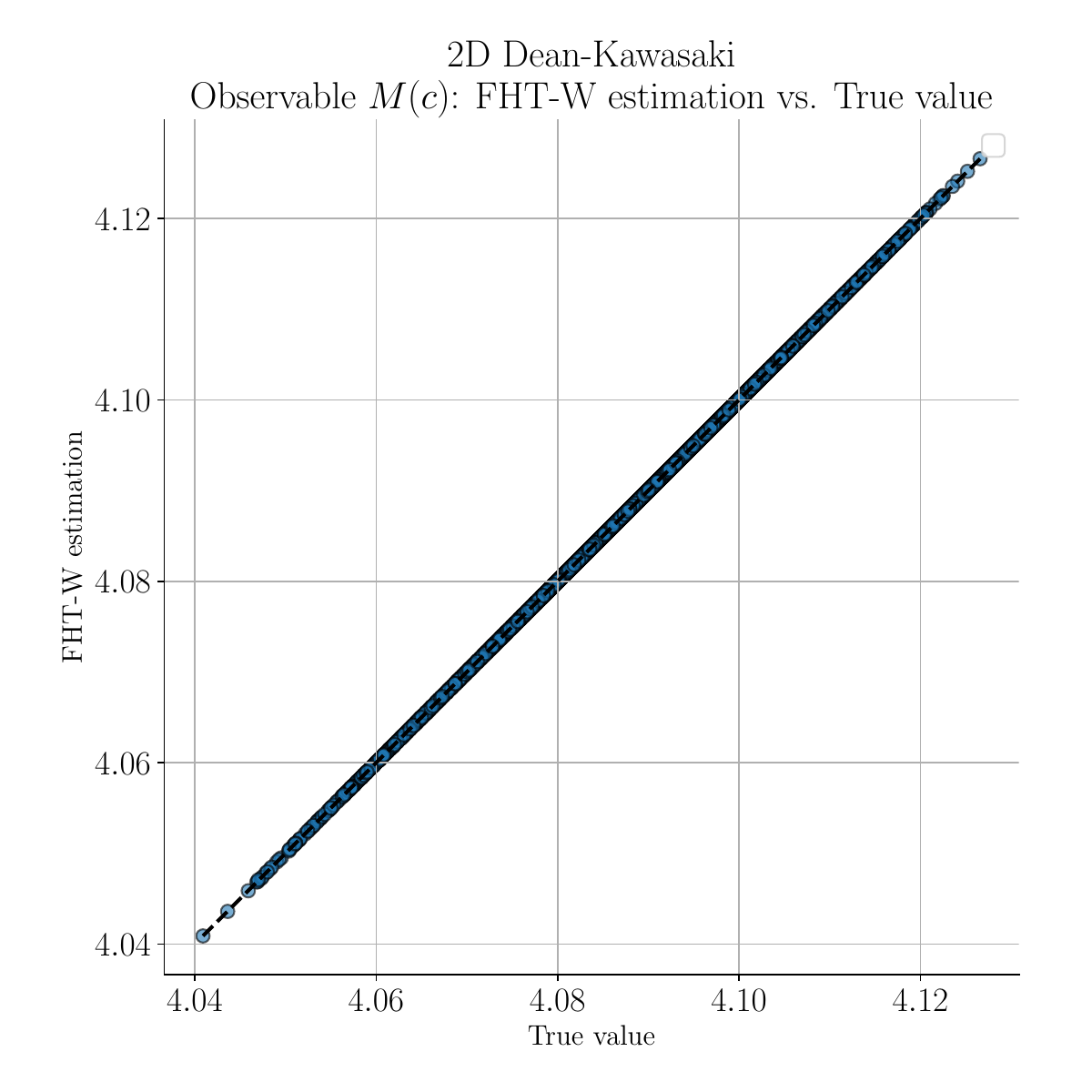}
    \caption{Interpolations for 1D Dean-Kawasaki model and 2D Dean-Kawasaki model without potential functions. The target function is \(M(c) =M(s^{-1}(\mathrm{IDWT}(c)))\), where \(M(\pi)\) is the Shannon entropy function. The plot compares the function evaluations of \(M(c)\) and the FHT-W ansatz obtained from the interpolation procedure in \Cref{alg:tree-based FTN interpolation internal}. The \(20{,}000\) evaluation points are chosen uniformly within the hypercube that contains the support of the data in the $s$ variable. The mean relative prediction error is \(8.7\times 10^{-9}\) for 1D Dean-Kawasaki and \(1.7\times 10^{-6}\) for 2D Dean-Kawasaki.}
    \label{fig:1D_DK_interpolation}
\end{figure}

We first use the FHT-W density estimate to predict the Shannon entropy
\[
M(\pi) = -\sum_{j = 1}^{d}\pi_j \log(\pi_j).
\]
As described in \Cref{sec: interpolation}, observable estimation requires constructing the composed function \(M(c) = M(\pi(c))\) via \Cref{alg:tree-based FTN interpolation internal}, using the interpolation parameters listed in \Cref{tab:parameters}. The result is shown in \Cref{fig:1D_DK_interpolation}, which confirms that the interpolation of \(M(c)\) is highly accurate.

The predicted entropy obtained after interpolation is \(4.15827\), while the true reference entropy from the Monte Carlo estimation is \(4.15825\). The relative FHT-W prediction error on \(M(\pi)\) is \(5.3\times 10^{-6}\). We repeat the experiments for the 2-R\'enyi entropy \(M(\pi) = -\log(\sum_{j= 1}^{d}\pi_j^2),\) and we likewise obtain a relative prediction error of \(1.0\times 10^{-5}\). For both observables, the FHT-W prediction error is on the same order as the Monte Carlo sampling error from $B = 6000$ trajectories.

\edit{The two-point correlation function is a standard observable for benchmarking density estimation methods on lattice models. In quantum state tomography, for instance, all leading methods are evaluated by their ability to capture the two-point correlation accurately (c.f. \cite{carrasquilla2019reconstructing, huang2020predicting}). Predicting the full correlation matrix from the estimated density is a stringent test because it requires the ansatz to represent the joint statistics of all pairs of sites simultaneously, rather than just marginals or scalar summaries. We therefore adopt the correlation matrix as a benchmark for the FHT-W density estimate.
The correlation matrix is defined by}
\[
R(i, i') =\mathrm{Corr}_{\Pi \sim P_{\pi}(\pi, t)}\left(\Pi_i, \Pi_{i'}\right),
\]
and we carry out the estimation by performing \(O(d^2)\) observable estimation tasks. \Cref{Fig: 1D covariance} shows that the predicted correlation closely matches the ground truth obtained from sample estimation.

\paragraph{Dean-Kawasaki with external potential and repulsion}
We consider the Dean-Kawasaki model with an external potential term \(V_1\) and a pairwise potential term \(V_2\). 
In this case, the equation reads
\[
\partial_t \pi = \frac{1}{\beta}\Delta \pi + \mathrm{div}\left(\sqrt{\frac{2\pi}{\beta N}}\eta + \pi\nabla V_1  + \pi\nabla V_2 * \pi\right),
\]
with parameters as in \Cref{tab:parameters}. For the external potential, we take
\[
V_1(x) = -15\cos(2\pi(x - 0.5)),
\]
which leads the distribution \(\pi(x, t)\) to be more concentrated around \(x \approx 0.5\). 

For the pairwise potential, we take \(V_2\) to be a regularized repulsion term written as follows:
\[
V_2(x) = -\frac{6}{x^2 + 0.01},
\]
which is a smooth penalization term. Due to the periodic boundary conditions, for two locations \(z, z' \in [0, 1]\), we choose their relative distance to be \(x = \min(|z-z'|, 1 - |z-z'|)\), and \(V_2(x)\) is the two-body interaction potential. We remark that \(V_2(x)\) and \(V_2'(x)\) are insignificant unless \(x\) is close to \(0\), and so boundary effects from this distance definition are negligible. 

The density estimation follows the same procedure as the zero-potential case, with parameters given in \Cref{tab:parameters}. The smaller time step ($d t = 0.0002$) is necessary due to the stiffer dynamics introduced by the potential terms.
We repeat the entropy estimation procedure. The FHT-W interpolation of \(M(\pi)\) in the \(c\) coordinate again achieves high precision: the relative prediction error is \(1.3\times 10^{-5}\) for the Shannon entropy and \(1.4\times 10^{-5}\) for the 2-R\'enyi entropy. Both errors are on the same order as the Monte Carlo sampling error from $B = 6000$ trajectories.

Finally, we use the FHT-W model to predict the correlation matrix. \Cref{Fig: 1D covariance Dean potential} shows that the FHT-W ansatz quite accurately predicts the correlation matrix \(R\).

\subsection{2D Dean-Kawasaki models}\label{sec: 2D numerical experiment}

\paragraph{Dean-Kawasaki without potential function}
We consider the 2D Dean-Kawasaki model with \(V_1 = V_2 = 0\),
\[
\partial_t \pi = \frac{1}{\beta}\Delta \pi + \mathrm{div}(\sqrt{\frac{2\pi}{\beta N}}\eta),
\]
with parameters as in \Cref{tab:parameters}. The initial condition is \(\pi(x, 0) = 1\) for all \( x \in [0, 1]^2\). We use a Legendre polynomial degree of $q = 15$ (rather than $q = 25$ as in 1D), as we empirically observe that the wavelet coefficients in the 2D model have a simpler univariate structure at each node.

We first test the observable estimation in terms of the Shannon entropy and the 2-R\'enyi entropy. The 2D wavelet transform leads to a more complicated functional form for \(M(c) = M(s^{-1}(\mathrm{IDWT}(c)))\). We show the result of FHT-W interpolation in \Cref{fig:1D_DK_interpolation}, which demonstrates that the interpolation of \(M(c)\) is highly accurate. For the Shannon entropy, the relative prediction error is \(5.5\times 10^{-6}\). For the 2-R\'enyi entropy, the relative prediction error is \(4.5\times 10^{-5}\).

\begin{figure}[h]
  \centering
  \includegraphics[width = \linewidth]{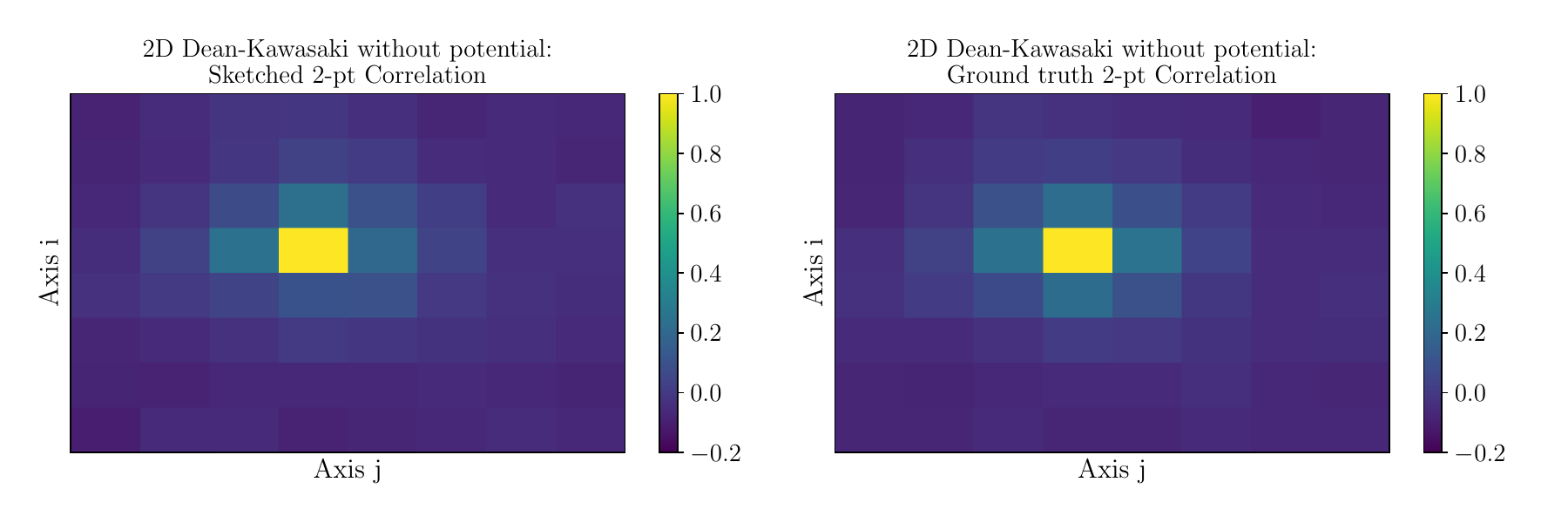}
  \caption{2D Dean-Kawasaki model with \(V_1 = V_2 = 0 \). The plot of the two-point correlation function predicted by the FHT-W ansatz. The mean prediction error is \(0.012\) and the max prediction error is \(0.128\).}
  \label{Fig: 2D_DK_2_pt_correlation_64}
    \includegraphics[width = \linewidth]{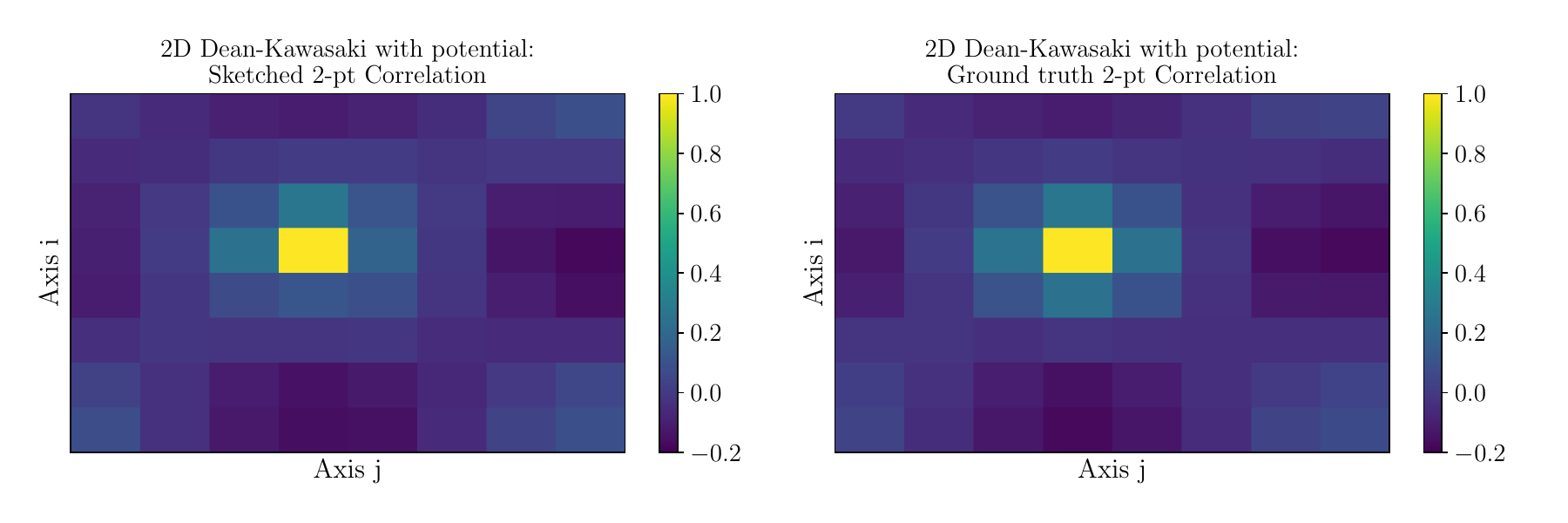}
  \caption{2D Dean-Kawasaki model with an external potential \(V_1\) and a pairwise potential term \(V_2\) modeling particle repulsion. The plot of the two-point correlation function predicted by the FHT-W ansatz. The mean prediction error is \(0.016\) and the max prediction error is \(0.119\).}
  \label{Fig: 2D_Dean_potential_2_pt_correlation_64}
\end{figure}

Lastly, we use the obtained FHT-W ansatz to predict the value of the two-point correlation function,
\[
f(i, j) =\mathrm{Corr}_{\Pi \sim P_{\pi}(\pi, t)}\left(\Pi_{(i,j)}, \Pi_{(4, 4)}\right),
\]
which we compute via repeated observable estimations. \Cref{Fig: 2D_DK_2_pt_correlation_64} shows that the proposed FHT-W ansatz closely matches the ground truth. 

\paragraph{Dean-Kawasaki with external potential and repulsion}
We consider the 2D Dean-Kawasaki model with an external potential term \(V_1\) and a pairwise potential term \(V_2\):
\[
\partial_t \pi = \frac{1}{\beta}\Delta \pi + \mathrm{div}\left(\sqrt{\frac{2\pi}{\beta N}}\eta + \pi\nabla V_1  + \pi\nabla V_2 * \pi\right),
\]
with parameters as in \Cref{tab:parameters}. For the external potential and pairwise potential, we take 
\[
V_1(x_1, x_2) = -15\cos(2\pi(x_1 - 0.5)) -15\cos(2\pi(x_2 - 0.5)),
\]
which leads the distribution \(\pi(x, t)\) to be more concentrated around \(x \approx (0.5, 0.5)\). For the pairwise interaction, we take 
\[
V_2(x_1, x_2) = -\frac{800}{x_1^2 + 0.1} -\frac{800}{x_2^2 + 0.1},
\]
which produces substantial repulsion between nearby particles. In this case, the dimension is \(d = 64\). The density estimation follows the parameters given in \Cref{tab:parameters}.

For the Shannon entropy and 2-R\'enyi entropy, we respectively obtain a relative error of \(3.5\times 10^{-5}\) and \(1.3\times 10^{-4}\). Lastly, we use the obtained FHT-W ansatz to predict the value of the two-point correlation function. \Cref{Fig: 2D_Dean_potential_2_pt_correlation_64} shows that the proposed FHT-W ansatz closely matches the ground truth.

\section{Conclusion}\label{sec: conclusion}
This work presents a computational framework for solving the Fokker-Planck equation associated with discretized Dean-Kawasaki models. The approach combines three ingredients: (1) a finite-volume discretization of the Dean-Kawasaki SPDE yielding a high-dimensional SDE, (2) a logarithmic coordinate transformation followed by a wavelet decomposition that maps the probability simplex to a Euclidean domain amenable to tensor network methods, and (3) a hierarchical Tucker decomposition \cite{hackbusch2009new} in the wavelet coordinate system (FHT-W) for density estimation and observable computation.

This work focuses on the strong-diffusion regime where all density components remain bounded away from zero. A sensitivity analysis of the bond dimension $r$ and polynomial degree $q$ is provided in the supplementary material, confirming that the method does not require delicate parameter tuning. Extending the approach to weaker diffusion regimes, developing adaptive strategies for $r$ and $q$, and applying iterative refinement techniques to improve upon the single-pass density estimation are directions for future work.

\section*{Code and data availability}
The code to reproduce all experiments and figures in this paper is publicly available at \url{https://github.com/Xun-Tang123/FHT_for_deans_equation}. The repository includes simulation scripts, density estimation routines, plotting scripts, precomputed reference data, and a \texttt{reproduce.sh} script that regenerates all results. A Docker container is provided for exact reproducibility of the computational environment. The repository contains automated tests verifying numerical correctness.

\bibliographystyle{siamplain} 
\bibliography{references}

\end{document}

\maketitle

\section{Sensitivity analysis: bond dimension $r$ and polynomial degree $q$}
\label{sec:sensitivity}

This supplement reports a sensitivity check for two approximation parameters in the FHT-W density-estimation pipeline:
the hierarchical bond dimension~$r$ and the univariate Legendre polynomial degree~$q$. We remark that $r$ is to be understood as a maximally allowed internal bond dimension, as the sketching algorithm adaptively chooses the rank based on the singular value decomposition that was done during density estimation. However, $r$ serves an upper bound on the rank to be chosen during the SVD.

\paragraph{Setup}
We use the 1D Dean--Kawasaki model with periodic boundary conditions and no external or pairwise potentials ($V_1=V_2=0$),
discretized with $d=64$ grid cells.
Unless otherwise stated, we adopt the same simulation/fit protocol as in the main manuscript:
$\beta=0.05$, $N=1000$, $d t=0.005$, $T=1$, and $B=6000$ trajectories.
We assess accuracy at $t=0.5$ using the two-point correlation (correlation matrix) of the finite-volume state $\Pi(t)$.

\paragraph{Error metrics (correlation prediction)}
Let $R\in\mathbb{R}^{d\times d}$ be a Monte Carlo reference correlation matrix computed from particle simulations,
and let $\widehat{R}$ be the correlation matrix predicted by the fitted FHT-W ansatz.
We report the max and mean absolute entrywise errors
\[
\mathrm{MaxErr} := \max_{1\le i,i'\le d}\bigl|\widehat{R}_{ii'}-R_{ii'}\bigr|,
\qquad
\mathrm{MeanErr} := \frac{1}{d^2}\sum_{i=1}^d\sum_{i'=1}^d \bigl|\widehat{R}_{ii'}-R_{ii'}\bigr|.
\]
These are the ``max prediction error'' and ``mean prediction error'' reported for correlation plots in the main manuscript.

\begin{table}[t]
\centering
\caption{Sensitivity of correlation prediction error at $t=0.5$.
Left: varying bond dimension $r$ with $q=25$ fixed. Right: varying polynomial degree $q$ with $r=15$ fixed.}
\label{tab:sensitivity_rq}
\begin{tabular}{@{}ccc|ccc@{}}
\toprule
\multicolumn{3}{c|}{Varying $r$ ($q=25$ fixed)} & \multicolumn{3}{c}{Varying $q$ ($r=15$ fixed)} \\
\cmidrule(r){1-3}\cmidrule(l){4-6}
$r$ & MaxErr & MeanErr & $q$ & MaxErr & MeanErr \\
\midrule
5  & 0.48 & 0.11 & 10 & 0.07 & 0.02 \\
10 & 0.04 & 0.02 & 15 & 0.05 & 0.02 \\
15 & 0.08 & 0.03 & 20 & 0.06 & 0.03 \\
20 & 0.06 & 0.01 & 25 & 0.09 & 0.03 \\
\multicolumn{3}{c|}{} & 30 & 0.14 & 0.06 \\
\bottomrule
\end{tabular}
\end{table}

\FloatBarrier

\paragraph{Takeaway}
In this test, the correlation errors are small for $r\ge 10$ and for $q$ in the range $10$--$25$.
The 1D density-fit parameters used in the main manuscript ($r=20$ and $q=25$) lie in this stable regime.